\documentclass[11pt]{article}
\usepackage{fullpage}
\usepackage{amsfonts}
\usepackage{amsthm}
\usepackage{graphicx}
\usepackage{epstopdf}
\usepackage{algorithmic}
\usepackage{wrapfig}
\usepackage{ragged2e}
\usepackage{tabularx}
\usepackage{mathdots}

\usepackage{soul}
\usepackage[normalem]{ulem}
\newcommand\redout{\bgroup\markoverwith
{\textcolor{red}{\rule[0.5ex]{2pt}{0.8pt}}}\ULon}

\ifpdf
  \DeclareGraphicsExtensions{.eps,.pdf,.png,.jpg}
\else
  \DeclareGraphicsExtensions{.eps}
\fi

\usepackage{mathtools}

\usepackage{xcolor}

\newcommand\bb{\mathbb}

\newcommand\R{\bb{R}}

\newcommand\p{\lambda}

\usepackage{bm}
\newcommand\xx{x}
\newcommand\ff{f}
\newcommand\FF{F}
\newcommand\vvarphi{\varphi}
\newcommand\g{g}
\newcommand\rr{r}
\newcommand\kkappa{\kappa}
\newcommand\FK{\Psi}

\newcommand{\flow}[2]{\vvarphi_{#1}(#2)}

\usepackage{enumitem}
\setlist[enumerate]{leftmargin=.5in}
\setlist[itemize]{leftmargin=.5in}

\theoremstyle{definition}
\newtheorem{definition}{Definition}
\newtheorem{proposition}{Proposition}
\newtheorem{theorem}{Theorem}

\usepackage{authblk}

\title{Continuation of fixed points and bifurcations\\ from ODE to flow-kick disturbance models\thanks{This work was supported by Carleton College's Towsley Endowment and New World Endowed Fund.}}

\author[*]{Katherine Meyer}
\affil[*]{Carleton College, Northfield, MN}
\author[*]{Horace Fusco} 
\author[*]{Collin Smith}
\author[**]{Alanna Hoyer-Leitzel}
\affil[**]{Mount Holyoke College}

\usepackage{amsopn}

\begin{document}

\maketitle

\vspace{0.1cm}
\begin{abstract}
    Ecosystems are shaped by disturbances such as fires and harvests that occur on faster timescales than recovery processes like population growth. Flow-kick models resolve these disturbances as discrete, periodic impulses in time, while ODE models smooth out disturbances in time as continuous processes. 
    Here we assess this simplification by comparing  dynamics of continuously disturbed ODE models to those of flow-kick models with the same average disturbance rate.
    In the case that kicks are small and high-frequency, we find that an ODE predicts similar dynamics as an analogous discrete disturbance model.  First, we prove that flow-kick maps generate an analogous vector field in the limit as the kick magnitude and return period approach zero in a fixed proportion. Second, we present conditions under which equilibria, saddle-node bifurcations, and transcritical bifurcations continue from ODE to analogous flow-kick systems. 
    We also provide numerical evidence that similarities between continuous and discrete disturbance models can break down as  kicks become less frequent and larger. We illustrate implications of these differences 
    in a nonspatial Klausmeier model of vegetation and precipitation dynamics. 
    We conclude that although ODEs may suffice to model high-frequency disturbances, resolving lower-frequency disturbances in time may be essential to effectively predicting their effects. 
\end{abstract}

%

\section{Introduction}
Modeling disturbances in an ecosystem---or any system---requires choices along a number of axes. Disturbances can be treated as deterministic or stochastic; spatially homogeneous or heterogeneous; temporally continuous or discrete. 
While a model of wildfires on the scale of 1,000 km might include stochastic components and spatial heterogeneity, controlled burns on the scale of 100 meters might be reasonably modeled as deterministic and spatially uniform. 
We focus on the case of deterministic, nonspatial models and probe the third axis: whether disturbances occur continuously in time or at discrete time points. 
In particular, we ask how this  choice impacts model predictions. Given that simple differential equations are used to model real-world disturbances that are markedly discrete in flavor, such as fire \cite{accatino2010tree} and rainfall \cite{klausmeier1999regular} events, do the predictions of discrete disturbance models differ enough to justify their increased complexity cost?

To explore this question, we model disturbances continuously with ordinary differential equations (ODEs) and discretely with impulsive differential equations \cite{lakshmikantham1989theory}. We focus on a subset of impulsive differential equations called flow-kick systems, which have emerged as a tool to study the interplay between discrete disturbance and continuous recovery when disturbances occur periodically in time \cite{meyer2018quantifying,zeeman2018resilience,hoyer2021impulsive}. A flow-kick map is built from two processes: a flow phase governed by an ODE and a discrete kick that delivers disturbance. This modeling structure appears under a variety of names in applications, including hurricane impacts on coral reefs \cite{ippolito2016alternative}, neurobiology of drug addiction \cite{chou2022mathematical}, fires \cite{tamen2016tree,tamen2017minimalistic,hoyer2021impulsive}, seasonal grazing \cite{ritchie2020episodic}, biocontrol \cite{fulk2022exploring}, and viral exposure \cite{hoyer2023immuno}. Theoretical studies have also uncovered rich dynamic structures such as horseshoes and strange attractors that emerge from periodically kicking simple flows \cite{wang2003strange,lin2010dynamics,naudot2021complexity}. Packaging one cycle of flow-kick into a map on state space opens up study of disturbance dynamics to the theory of iterated maps, including fixed points (where disturbance and recovery balance) and bifurcations (potential tipping points). Unlike the standard study of maps, however, one typically lacks an exact closed form for the flow-kick map. 

This added complexity of the discrete approach both motivates our comparison to continuous systems and poses a challenge for analyzing flow-kick maps---a challenge which we address analytically using Taylor expansions of the flow and numerically using the MatcontM continuation package. Using these tools, we compare continuous and discrete systems with the same average disturbance rate. We find that predictions are similar between continuous and discrete, high-frequency disturbance models, but that this similarity breaks down for lower frequencies. Thus using a discrete disturbance model may yield novel insights.\\

We outline our five main analytic results (boldfaced below) in the simple context of the logistic model of population growth 
\begin{equation}\label{eq:plainlog}
    x'=x(1-x)
\end{equation}
disturbed by harvesting. 

Figure 1 illustrates three disturbed logistic systems that feature the same average disturbance rate over time. 

\begin{figure}[b!]
    \centering \includegraphics[width=0.9\textwidth]{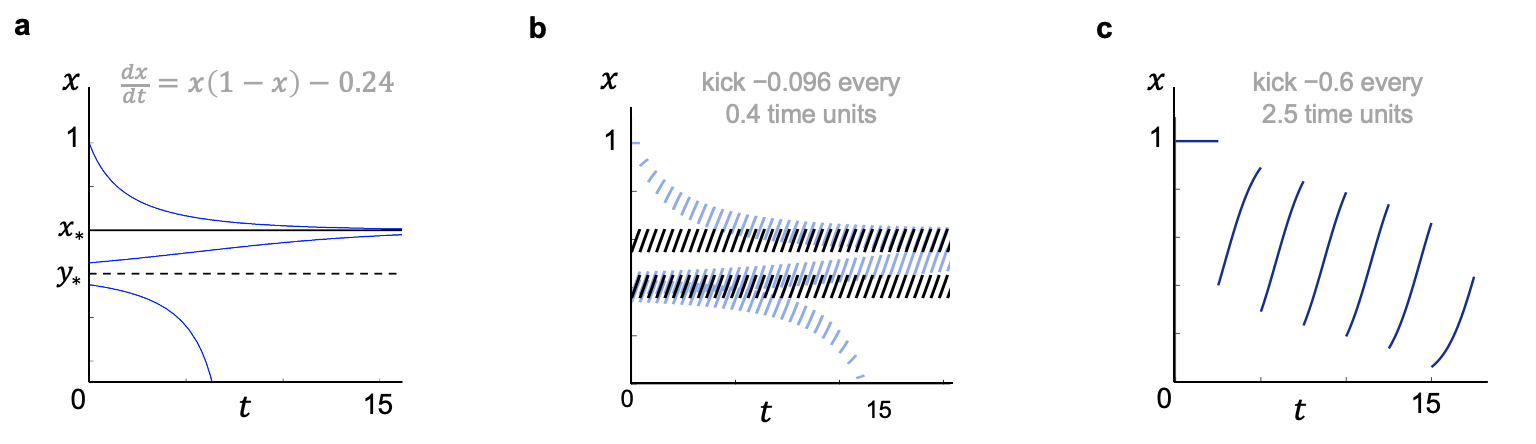}
        \caption{\small (a) Flow of the logistic ODE minus continuous   harvesting rate of $0.24$. (b, c) Flows according to undisturbed logistic equation with kicks delivered at the same average disturbance rate: $\text{kick}/(\text{flow time})=-0.24$.}
    \label{fig:logistic}
\end{figure}

In panel (a), the ODE
\begin{equation}\label{eq:distlog}
    x'=x(1-x)-0.24
\end{equation}
represents continuous disturbance at a rate of $-0.24$, which yields a stable equilibrium $x_*$ and an unstable equilibrium $y_*$. In panel (b), undisturbed logistic flows are interrupted every 0.4 time units by a kick of $-0.096$, chosen to maintain the average disturbance rate
\begin{equation}
    \frac{\text{kick}}{\text{flow time}}=\frac{-0.096}{0.4}=-0.24.
\end{equation} 
\noindent The continuous and high-frequency discrete disturbance models in Figure \ref{fig:logistic}ab exhibit similar dynamics: in particular, the flow-kick map illustrated in Figure \ref{fig:logistic}b has fixed points near $x_*$ and $y_*$ with corresponding stabilities. In Section \ref{sec:fixedpts} we prove that this phenomenon---\textbf{fixed point continuation from ODEs to flow-kick maps}---occurs under quite general conditions for $x\in\R^n$ (Proposition \ref{prop:fixedpt_n}). The disturbance rate (here a constant $-0.24$) can in fact be any smooth function of $x$. Furthermore, \textbf{fixed points inherit stability} under slightly stronger assumptions (Theorem \ref{thm:stability_n}). The proofs of Proposition \ref{prop:fixedpt_n} and Theorem \ref{thm:stability_n} use the Implicit Function Theorem, so we are guaranteed continuation only for sufficiently small flow times. 

Indeed, Figure \ref{fig:logistic}c shows that longer flow times of $2.5$ between kicks of $-0.6$ destroy the stable and unstable fixed points observed in the continuous and high-frequency disturbance models, despite having the same average disturbance rate. As noted in \cite{meyer2018quantifying,zeeman2018resilience} in the context of fishery management, infrequent harvesting at rates that appear sustainable in a continuous model could in fact crash the population. This difference underscores the importance of resolving discrete disturbances in time. \\

\begin{figure}[h!]
    \centering
\includegraphics[width=0.9\textwidth]{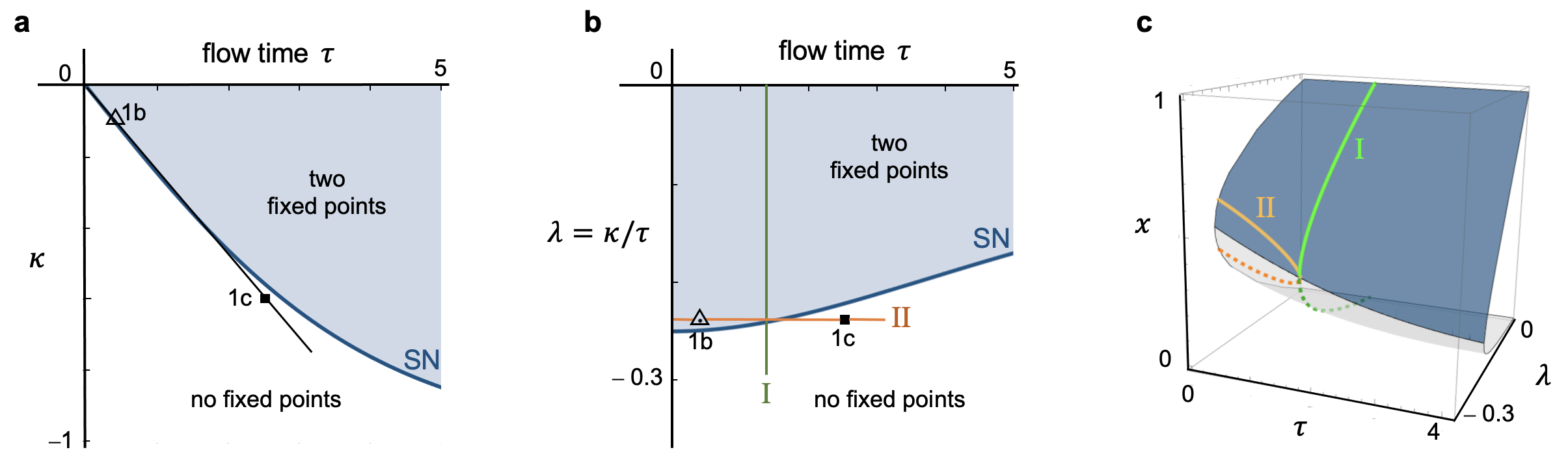}
    \caption{\small Stability diagram for the logistic flow-kick system in parameters flow time $\tau$ and kick $\kappa$ (panel a) or $\tau$ and disturbance rate $\p=\kappa/\tau$ (panel b). Marked points correspond to parameter combinations in Figure 1b and 1c. (c) Surface of flow-kick fixed points over $\tau,\p$ parameter space. The flow-kick map exhibits saddle-node bifurcations in both $\p$ (I) and $\tau$ (II).}
    \label{fig:logbifnew}
\end{figure}

The disappearance of the stable and unstable flow-kick fixed points from panels (b) to (c) in Figure \ref{fig:logistic} suggests a saddle-node bifurcation, reminiscent of the saddle-node bifurcation that occurs in the continuously disturbed system 
\begin{equation}\label{eq:logctsdis}
    x'=x(1-x)+\p
\end{equation}
at $\p=-0.25$. Because the logistic differential equation is analytically tractable, we can compute flow-kick fixed points and bifurcation curves exactly. Figure \ref{fig:logbifnew} presents a curve of saddle-node bifurcations in two different parameter spaces. 

In Figure \ref{fig:logbifnew}a we take the flow time $\tau$ and kick $\kappa$ as parameters for the flow-kick map, as used in previous papers \cite{meyer2018quantifying, hoyer2021impulsive, zeeman2018resilience}. 
As expected, the disturbance patterns corresponding to Figures \ref{fig:logistic}b and \ref{fig:logistic}c fall on opposite sides of the bifurcation curve.
We change coordinates in Figure \ref{fig:logbifnew}b,c to the parameters $\tau$ and $\p\equiv\kappa/\tau$, and again the patterns corresponding to Figures \ref{fig:logistic}b and \ref{fig:logistic}c fall on opposites of the bifurcation curve.
Note that families of flow-kick maps with the same average disturbance rate $\kappa/\tau$ fall on lines through the origin in Figure \ref{fig:logbifnew}a and on horizontal lines in Figure \ref{fig:logbifnew}b. 

We can use either Figure \ref{fig:logbifnew}a or Figure \ref{fig:logbifnew}b to visualize the limit as we move from the flow-kick system of Figure 1b towards the ODE of Figure 1a, taking the flow time and kick to zero in a fixed proportion. This limit appears in Figure \ref{fig:logbifnew}a as movement towards the origin along a line with slope $-0.24$, and in Figure \ref{fig:logbifnew}b as movement towards the vertical axis along a horizontal line $\lambda=-0.24.$
We claim that in the limit
as $\tau$ and $\kappa$ approach zero in a fixed ratio $\kappa/\tau=\lambda_0$,
flow-kick systems 
approach the continuous system 
\begin{equation}\label{eq:ctslimitode}
    x'=x(1-x)+\p_0.
\end{equation} 
We call \eqref{eq:ctslimitode} the \textbf{continuous limit of the collection of flow-kick maps} with average disturbance rate $\p_0$. Proposition \ref{prop:ctslimit} makes this idea precise and generalizes it to a broader class of  disturbance functions.

Figures \ref{fig:logbifnew}a and \ref{fig:logbifnew}b thus contain an essential nuance at the origin and vertical axis, respectively. At the origin in Figure \ref{fig:logbifnew}a, $\tau=\kappa=0$ and the flow-kick map is the identity. However, the continuous limit result in Proposition \ref{prop:ctslimit} suggests we might imagine a continuum ODE systems \eqref{eq:ctslimitode} packed singularly at the origin in Figure \ref{fig:logbifnew}a---one for each value of $\lambda_0$. In Figure \ref{fig:logbifnew}b we can imagine this continuum of ODEs spread along the vertical $\lambda$ axis.
The $(\tau,\lambda)$ coordinates used in Figure \ref{fig:logbifnew}b provide not only the visual advantage of unpacking distinct limiting ODEs but also an analytic tool for the continuation proofs in this paper. Definition \ref{def:fkr} presents these coordinates for flow-kick maps of a more general form.

Figure \ref{fig:logbifnew}b also highlights relationships between bifurcations in the continuously and discretely disturbed logistic system.
Given that each point $(\tau,\p)=(0,\p_0)$ on the vertical axis in Figure \ref{fig:logbifnew}b represents---in some sense---an ODE of the form \eqref{eq:ctslimitode}, we observe that the curve of flow-kick saddle-node bifurcations (marked SN) emanates from the saddle-node bifurcation at $\p=-0.25$ for the continuously disturbed logistic system. We prove in Section \ref{sec:bif} that for $x\in\R^n$, \textbf{saddle-node bifurcations continue from ODE to flow-kick systems} (Theorem \ref{thm:s-n}). We generalize from the constant disturbance rate $\p$ considered here to a parameterized disturbance rate function $r(x,\p)$. Theorem \ref{thm:s-n} connects bifurcations in $\p$ of a continuous system to bifurcations in $\p$ of a flow-kick system, shown for the logistic system as the green curve (I) in Figure \ref{fig:logbifnew}c. 

Transcritical bifurcations also commonly arise in ecological models that feature $x=0$ as an extinction equilibrium that is unchanged by population growth or disturbance processes. In Section \ref{sec:bif}, we show that \textbf{transcritical bifurcations continue} from ODE to flow-kick systems in this setting (Theorem \ref{thm:tc}). In sum, our local continuation results suggest that ODEs have a good chance of capturing some of a flow-kick system's dynamic features when discrete disturbances  are  frequent and small. 

The paper is organized as follows. Section \ref{sec:prelim} establishes definitions and notation. In section \ref{sec:ctsanalog} we formalize the idea of a continuous limit of flow-kick maps that maintain a constant average disturbance rate. Sections \ref{sec:fixedpts} and \ref{sec:bif} treat continuation of fixed points and bifurcations, respectively, from ODEs to analogous flow-kick models of disturbance. In section \ref{sec:examples} we illustrate our results numerically in two models pertaining to vegetation-water and predator-prey dynamics. These examples  highlight the possibility of qualitatively different behavior between ODE and analogous flow-kick models of disturbance when disturbances are large and low-frequency. They also raise topics for future study such as continuation of Hopf bifurcations, which we discuss in Section \ref{sec:disc}.

\section{Preliminaries and Notation}\label{sec:prelim}
\subsection*{Flows}
Throughout, we consider undisturbed dynamics modeled by an ODE of the form
\begin{equation}\label{eq:undist}
    x'=\ff(\xx)
\end{equation}
where $\xx(t)$ is a function from $\R$ to $\R^n$, $'$ denotes $d/dt$, and $\ff:U\to\R^n$ is a $C^4$ vector field defined on an open set $U\subseteq \R^n$. 
The vector field $\ff$ generates a local $C^4$ flow function on a subset of $\R\times U$ given by 
\begin{equation}   
    \vvarphi(t,\xx_0)=\xx(t),
\end{equation}
where $x(t)$ is the position of a solution to \eqref{eq:undist} that starts at $\xx_0$ and flows for time $t$. We will assume that the flow is defined on any time interval of interest. Fixing $t$ yields a time-$t$ map $\vvarphi_t:U\to U$ given by
\begin{equation}\label{eq:timetmap}
\vvarphi_t(\xx)=\vvarphi(t,\xx).
\end{equation}
Given an ODE $x'=g(x)$, we will use the terms equilibrium and fixed point interchangeably to describe $x\in\R^n$ such that $g(x)=0$.

\subsection*{Kicks}
To incorporate disturbances into dynamics, we choose a flow time $t=\tau$, representing the period between disturbances, and add a discrete kick to the time-$\tau$ map. The following definition generalizes several found in the literature (e.g. \cite{hoyer2021impulsive,meyer2018quantifying,zeeman2018resilience}).
\begin{definition}\label{def:fk}
    A flow-kick map $\FK_{\tau,k}:U\to \R^n$ with kick function $k: U \to \R^n$ is given by  \[\FK_{\tau,k}(\xx)=\vvarphi_\tau(\xx)+k(\vvarphi_\tau(\xx))\]
    where $\varphi_{\tau}$ is the time-$\tau$ map \eqref{eq:timetmap} corresponding to the ODE \eqref{eq:undist} and the kick function 
    $k$ is $C^2$.
\end{definition}


\noindent For example, if $k(u)=\kkappa$ for a constant $\kkappa\in\R^n$, then the flow-kick map
\begin{equation}
\FK_{\tau,k}(\xx)=\vvarphi_\tau(\xx)+\kkappa
\end{equation}
matches the disturbance forms considered in \cite{meyer2018quantifying}. By choosing $k(u)=-\kappa u$ we obtain the ``multiplicative" flow-kick map studied in \cite{zeeman2018resilience}:
\begin{equation}
    \FK_{\tau,k}(\xx)=\vvarphi_\tau(\xx)-\kappa\vvarphi_\tau(\xx)
\end{equation}

\subsection*{Discrete Dynamics}

Iterating $\FK_{\tau,k}$ generates discrete dynamics. We focus the present study on the fixed points of these dynamics, where disturbance and recovery processes balance. 

\begin{definition}
    A flow-kick fixed point is a solution to $\FK_{\tau,k}(\xx)=\xx$, or equivalently \\$\vvarphi_{\tau}(\xx)+k(\varphi_\tau(\xx))=\xx$.
\end{definition}


It is worth noting that since $\FK_{\tau,k}$ maps from $U$ to $\R^n$, it is possible for a flow-kick trajectory to leave $U$. For example, if $U$ is the interval $(0,\infty)$ in the introductory example of harvested logistic growth, the bottom trajectory in Figure \ref{fig:logistic}b leaves $U$. However, this will not trouble our local study of fixed points in the open set $U$.

\subsection*{Disturbance Coordinates}
In the context of the logistic system with constant kicks $\kappa$, Figure \ref{fig:logbifnew} summarizes system behavior over two alternative  disturbance parameter choices: $(\tau,\kappa)$ in panel (a) and $(\tau, \p\equiv\kappa/\tau)$ in panel (b).
We generalize the coordinate change between these two panels to non-constant disturbances. A disturbance rate function can be incorporated continuously into an ODE model as $x'=f(x)+r(x)$ or used to generate a kick 
$k(u)=\tau r(u)$
in a flow-kick model. In modeling a disturbance pattern, kick magnitude $k(u)$ and flow time $\tau$ can be chosen independently.
The disturbance rate function $r(u)$ gives the average disturbance rate over time and facilitates comparison of a discrete disturbance model and its analogous continuous disturbance model. 

\begin{definition}\label{def:fkr}
    The flow-kick map associated with a disturbance rate function $r:U\to\mathbb{R}^n$ and flow time 
    $\tau$ is \[\Phi_{\tau,r}(x)=\FK_{\tau,\tau r}=\vvarphi_\tau(x)+\tau r(\varphi_\tau(x)).\]
    where $\varphi_{\tau}$ is the time-$\tau$ map \eqref{eq:timetmap} corresponding to the ODE \eqref{eq:undist} and the disturbance rate function 
    $r$ is $C^2$.
\end{definition}

\noindent Henceforth we use Definition \ref{def:fkr} to describe families of flow-kick maps in terms of their disturbance rate functions $r$ and flow times $\tau$. The disturbance rate function $r$ is omitted from the subscript of $\Phi$ when it is clear from context. The results that follow apply to  $\Phi_{\tau,r}$ (Definition \ref{def:fkr}) rather than $\Psi_{\tau,k}$ (Definition \ref{def:fk}).

\section{Continuous limit of flow-kick systems}\label{sec:ctsanalog}

The similarities between the continuous harvest model (Figure \ref{fig:logistic}a) and the flow-kick harvest model (Figure \ref{fig:logistic}b)  reflect a broader phenomenon: flow-kick maps generate an analogous vector field in the limit as flow times go to zero and a constant disturbance rate is maintained.
This result, which has been shown for specific flow-kick maps in \cite{hoyer2021impulsive, meyer2018quantifying}, provides a bridge between flow-kick models and purely ODE models with the same average disturbance rate, a common quantity used in ecological modeling. Here we state and prove a more general formulation using the current notation.

\begin{proposition}\label{prop:ctslimit}
    Let $f(\xx)$ be a vector field governing undisturbed dynamics and $r(\xx)$ be a disturbance rate function. In the limit as $\tau\to 0$, the family of flow-kick maps $\Phi_{\tau,r}(x)=\vvarphi_\tau(\xx)+\tau r(\vvarphi_\tau(\xx))$ generates the vector field $f(x)+r(x)$.
\end{proposition}

\begin{proof}
To calculate the vector field generated by $\Phi_{\tau,r}$ as $\tau\to 0$ we use the derivative of this time-dependent mapping with respect to $\tau$ at $\tau=0$. Note that $\Phi_{0,r}(\xx)$ is the identity; that is, $\Phi_{0,r}(\xx)=x$.  We have
\begin{align}
    \lim\limits_{\tau\to 0}\frac{\Phi_{\tau,r}(\xx)-x}{\tau}
    &=
        \lim\limits_{\tau\to 0}\frac{\vvarphi_\tau(\xx)+\tau r(\vvarphi_\tau(\xx))-x}{\tau}.
       \label{eq:diffquot}   
\end{align}
\noindent Because the flow function $\vvarphi$ is twice differentiable, it can be expanded in the time variable via Taylor's formula as
\begin{equation}\label{eq:taylor}
\begin{aligned}
    \vvarphi(\tau,\xx)&= \vvarphi(0,\xx)+\tau \frac{\partial\vvarphi}{\partial t}(0,x)+O(\tau^2)\\
    &=\xx+\tau f(x)+ O(\tau^2).
\end{aligned}
\end{equation}
Here we have employed Landau big-O notation in a neighborhood of $\tau=0$. 
Substituting \eqref{eq:taylor} into \eqref{eq:diffquot} and simplifying yields
\begin{equation}
    \begin{aligned}
    \lim\limits_{\tau\to 0}\frac{\Phi_{\tau,r}(\xx)-x}{\tau}
    &=\lim\limits_{\tau\to 0}\frac{\tau f(\xx)+\tau r(\vvarphi_\tau(\xx))+O(\tau^2)}{\tau}  \\
    &= \lim\limits_{\tau\to 0} \big(f(\xx)+ r(\vvarphi_\tau(\xx))+O(\tau)\big)\\
    &=f(\xx)+r(\xx),
    \end{aligned}
\end{equation}
as claimed. 
\end{proof}

We can visualize Proposition \ref{prop:ctslimit} for the harvested logistic growth example in Figure \ref{fig:logbifnew}, where the disturbance rate function is a constant $r(x)=\p$. Each value of $\p=\kappa/\tau$ yields a line with slope $\p$ through the origin in Figure \ref{fig:logbifnew}a. In Figure \ref{fig:logbifnew}b, fixing $\p=0.24$ yields a horizontal line. Proposition \ref{prop:ctslimit} concerns the limit of the flow-kick maps parameterized along these lines as $\tau\to0^+$. Note that limits along different lines generate different continuous systems $x'=f(x)+\p$. In this sense the origin in Figure \ref{fig:logbifnew}a and the vertical axis in Figure \ref{fig:logbifnew}b each represent a continuum of vector fields that result from taking $\tau$ to $0$ with different disturbance rates $\p$. On the other hand, the literal flow-kick map (Definition \ref{def:fkr}) parameterized by $\tau=0$ is the identity map, as $\Phi_{0,r}(\xx)=\vvarphi_0(x)+0r(\vvarphi_0(x))=x$. 

\section{Continuation of equilibria to fixed points}\label{sec:fixedpts}
Having established a connection between flow-kick models and an ODE models with the same average disturbance rate in Proposition \ref{prop:ctslimit}, we turn to the question of whether the ODE model captures dynamics of its flow-kick analogs. Proposition \ref{prop:fixedpt_n} and Theorem \ref{thm:stability_n} answer affirmatively in the case a hyperbolic fixed points and small, frequent kicks.
First, Proposition \ref{prop:fixedpt_n} shows that
under suitable nondegeneracy conditions, an equilibrium for continuous disturbance continues locally to a flow-kick fixed point when we discretize the disturbance, maintaining the same average rate:

\begin{proposition}\label{prop:fixedpt_n}
If $\xx_*\in U$ is an equilibrium for the continuous system $\xx'=\ff(\xx)+\rr(\xx)$ and $D_{\xx}[\ff+\rr](\xx_*)$ is invertible, then $x_*$ continues locally to flow-kick fixed points. That is, there exists $\delta>0$ and a family $\{\xx_\tau\}_{\tau\in[0,\delta)}\subset U$  satisfying  
\begin{enumerate}
    \item[(i)] each $\xx_\tau$ is a fixed point for the flow-kick map $\Phi_{\tau,r}(\xx)=\flow{\tau}{\xx}+\tau\rr(\flow{\tau}{\xx})$,
    \item[(ii)] the flow-kick fixed points $\xx_\tau$ vary continuously as a function of $\tau$, and
    \item[(iii)] the branch of flow-kick fixed points emerges from the equilibrium $\xx_*$; that is, $\xx_0=\xx_*$.
\end{enumerate}  
\end{proposition}

\begin{proof} Flow-kick fixed points $\xx\in U$ are those that satisfy $\Phi_{\tau,r}(\xx)=\xx$, or equivalently
\begin{equation}\label{eq:initial}
    \flow{\tau}{\xx}+\tau \rr(\flow{\tau}{\xx}) -\xx=0.
\end{equation}
When $\tau=0$, all points $\xx\in U$ satisfy equation (\ref{eq:initial}). To remove this degeneracy and connect to $\xx_*$ we expand $\vvarphi$ in $t$ about 0:  
\begin{align}
   \flow{\tau}{\xx}=\vvarphi(\tau,\xx)&=\varphi(0,\xx) + \tau \frac{\partial \vvarphi}{\partial t} (0,\xx) + \tau^2\int_0^1(1-u)\frac{\partial^2 \vvarphi}{\partial t^2}(u\tau, \xx) du \label{eq:taylormessy}\\
   &=\xx +\tau \ff(\xx) + \tau^2\int_0^1(1-u)\frac{\partial^2 \vvarphi}{\partial t^2}(u\tau, \xx) du. \label{eq:taylorclean}
\end{align}
In (\ref{eq:taylormessy}) we have employed an integral form of the remainder articulated by Folland \cite{folland1990remainder} (see Appendix \ref{sec:intform}). Substituting (\ref{eq:taylorclean}) into (\ref{eq:initial}) and simplifying yields 
\begin{equation}\label{eq:factors}
    \tau\left(\ff(\xx)+\rr(\flow{\tau}{\xx})+\tau\int_0^1(1-u)\frac{\partial^2 \vvarphi}{\partial t^2}(u\tau, \xx) du\right)=0.
\end{equation}
To solve (\ref{eq:initial}) it suffices to find zeros of the second factor in (\ref{eq:factors}). We do this via the Implicit Function Theorem applied to the function $\FF:\R\times U\to\R^n$ given by
\begin{equation}\label{eq:F}
    \FF(\tau,\xx)=\ff(\xx)+\rr(\flow{\tau}{\xx})+\tau\int_0^1(1-u)\frac{\partial^2 \vvarphi}{\partial t^2}(u\tau, \xx) du.
\end{equation}
Given the smoothness assumptions on $\ff$, $\rr$, and $\vvarphi$, one can confirm from (\ref{eq:F}) that $\FF$ is $C^1$. 
Further,  $\FF(0,\xx_*)=\ff(\xx_*)+\rr(\xx_*)=0$, since $\xx_*$ is a fixed point for the continuous system. 
Lastly, we have that 
\begin{align}
D_{\xx}\FF(0,\xx_*)&=D_{\xx}\left.[\ff(\xx)+\rr(\flow{0}{\xx})]\right|_{\xx=\xx_*}\\
&=D_{\xx}[\ff+\rr](\xx_*) 
\end{align}
is invertible, by hypothesis. It follows from the Implicit Function Theorem that for some $\delta>0$, there exists a $C^1$ function $\g:[0,\delta)\to U$ satisfying $\g(0)=\xx_*$ and $\FF(\tau,\g(\tau))=0$. Letting $\xx_\tau=\g(\tau)$, the proposition follows.
\end{proof}

With slightly stronger hypotheses on the spectrum of $D_{\xx}[\ff+\rr](\xx_*)$, the stability of continuous-disturbance fixed point $\xx_*$ persists locally along the branch $\xx_\tau$ of flow-kick fixed points:

\begin{theorem}\label{thm:stability_n}
If $\xx_*\in U$ is a hyperbolic equilibrium of the continuous system $\xx'=\ff(\xx)+\rr(\xx)$ and $\tau>0$ is sufficiently small, then flow-kick fixed points $\xx_\tau$ not only continue from $\xx_*$ as in \emph{Proposition \ref{prop:fixedpt_n}}  but also inherit the stability of $\xx_*$.
\end{theorem}

\begin{proof}
Consider a fixed point $\xx_*$ of $\xx'=\ff(\xx)+\rr(\xx)$ that is hyperbolic, meaning all eigenvalues of $D_{\xx}[\ff+\rr](\xx_*)$ have nonzero real part. By Proposition \ref{prop:fixedpt_n}, a continuous family of flow-kick fixed points $\{\xx_\tau\}_{\tau\in[0,\delta)}$ emanates locally in $\tau$ from $\xx_0=\xx_*$. We classify the stability of a fixed point $\xx_\tau$ using the Jacobian matrix
\begin{align*}
    D_{\xx}\Phi_{\tau,r}(\xx)&=D_{\xx}[\flow{\tau}{\xx}+\tau \rr(\flow{\tau}{\xx})]\\
    &=D_{\xx}\left[\xx+\tau\ff(\xx)+ \tau^2\int_0^1(1-u)\frac{\partial^2 \vvarphi}{\partial t^2}(u\tau, \xx) du\right]+D_{\xx}\Big[\tau\rr(\flow{\tau}{\xx})\Big]\\
    &=I + \tau  \left[D_{\xx}\ff(\xx)+D_{\xx}\rr(\flow{\tau}{\xx})+\tau D_{\xx} \int_0^1(1-u)\frac{\partial^2 \vvarphi}{\partial t^2}(u\tau, \xx) du \right]\\
    &=I+\tau A(\tau,\xx)
\end{align*}
where $A$ is the matrix-valued function
\[
A(\tau,\xx)=D_{\xx}\ff(\xx)+D_{\xx}\rr(\flow{\tau}{\xx})+\tau D_{\xx} \displaystyle\int_0^1(1-u)\frac{\partial^2 \vvarphi}{\partial t^2}(u\tau, \xx) du .
\]

Let $\ell_1(\tau),\dots,\ell_n(\tau)$ be the eigenvalues of $A(\tau,\xx_\tau)$, counted with multiplicities. Let $a_i(\tau)$ and $b_i(\tau)$ denote the real and imaginary parts of $\ell_i(\tau)$, respectively. The hyperbolicity of $\xx_*$ implies that each eigenvalue $\ell_i(0)$ of $A(0,\xx_0)=D_{\xx}\ff(\xx_*)+D_{\xx}\rr(\xx_*)$ has $a_i(0)\not=0$. 
Furthermore, $\ff$, $\rr$, $\vvarphi$, and $\xx_\tau$ are sufficiently smooth so that the function $A(\tau,\xx_\tau)$ varies continuously in $\tau$ along the branch of flow-kick fixed points.
It follows from continuity of eigenvalues with respect to matrix perturbations (see, e.g. \cite{stewart1990matrix}) that for some $\delta_1\in(0,\delta]$, $0\leq \tau<\delta_1$ ensures  $\text{sign}(a_i(\tau)) = \text{sign}( a_i(0))$ for $i=1,\dots,n$. In other words, $\delta_1$ bounds the neighborhood in which the eigenvalues of $A$ do not cross the imaginary axis. 

One may confirm via characteristic equations that the eigenvalues $m_1(\tau),\dots,m_n(\tau)$ of the derivative matrix $D_{\xx}\Phi_{\tau}(\xx_\tau)$ are related to the eigenvalues  of $A(\tau,\xx_\tau)$ by 
\begin{equation}\label{eq:eigs}
    m_i(\tau)=1+\tau \ell_i(\tau).
\end{equation}
Stability of the flow-kick fixed point $\xx_\tau$ depends on size of the moduli
\begin{equation}\label{eq:modulus}
|m_i(\tau)|=|1+\tau(a_i(\tau)+b_i(\tau)i)|
=\sqrt{1+\tau(2a_i(\tau)+\tau(a_i(\tau)^2+b_i(\tau)^2))}
\end{equation}
relative to 1.
Recall flow-kick fixed points correspond to $\tau>0$. In the case that $a_i(0)>0$, then for $0<\tau<\delta_1$, we have $a_i(\tau)>0$, and it follows readily that $|m_i(\tau)|>1$. 
If instead $a_i(0)<0$, then for $0<\tau<\delta_1$, we have $a_i(\tau)<0$. Note that the factor $2a_i(\tau)+\tau(a_i(\tau)^2+b_i(\tau)^2)$ in (\ref{eq:modulus}) is a continuous function of $\tau$ that evaluates to $2a_i(0)<0$ when $\tau=0$. Thus there exists a $\delta_2>0$ such that $2a_i(\tau)+\tau(a_i(\tau)^2+b_i(\tau)^2)<0$ whenever $0<\tau<\delta_2$. Taking $\delta_3=\min(\delta_1,\delta_2)$ ensures that $|m_i(\tau)|<1$ whenever $a_i(\tau)<0$ and $0<\tau<\delta_3$.

In this way stable and unstable eigenvalues $\ell_i(0)$ of the vector field $\ff(\xx)+\rr(\xx)$ at $\xx_*$ continue locally to stable and unstable, respectively, eigenvalues $m_i(\tau)$ of the flow-kick map $\Phi_{\tau,r}$ at $\xx_\tau$. 
\end{proof}

\phantom{.}\begin{wrapfigure}{r}{0.35\textwidth}
\vspace{-1cm}
\includegraphics[width=0.9\linewidth]{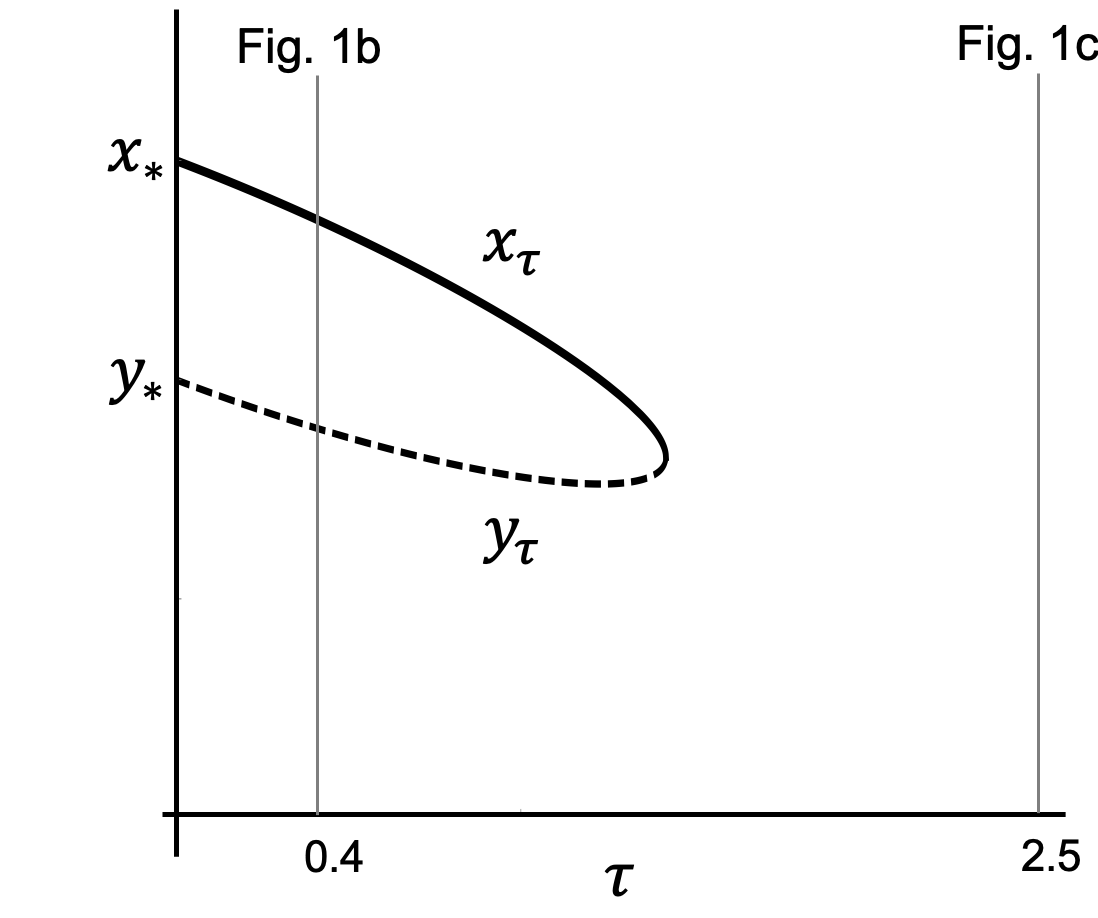} 
\caption{{\small{Logistic flow-kick fixed points with $r(u)=-0.24$}}}
\label{fig:cor1}
\vspace{1cm}
\end{wrapfigure}

\vspace{-1cm}
To connect the notation of Theorem \ref{thm:stability_n} to the introductory logistic example, consider logistic growth generated by \begin{equation*}
   x'=f(x)=x(1-x) 
\end{equation*} with disturbance rate $r(u)=-0.24$. 
Figure \ref{fig:cor1} illustrates the branches $x_\tau$ (stable) and $y_\tau$ (unstable) of flow-kick fixed points that continue from the ODE equilibria $x_*$ and $y_*$ shown in Figure \ref{fig:logistic}a. Vertical lines mark the values of $\tau$ used for flow-kick systems in Figure \ref{fig:logistic}b and \ref{fig:logistic}c.  
The continuous disturbance model predicts the qualitative dynamics of the frequent, small harvests ($\tau\!=\!0.4$, $\kappa\!=\!-0.096$) but not of the infrequent, large harvests ($\tau\!=\!2.5$, $\kappa\!=\!-0.6$) that occur at the same average rate. This is because
near $\tau=1.4$ the branches of fixed points meet in a saddle-node bifurcation---the same one illustrated by curve (II) in Figure \ref{fig:logbifnew}c.
We turn to bifurcation continuation in Section \ref{sec:bif}. 

\section{Continuation of bifurcations}\label{sec:bif}
As bifurcations in flow-kick systems represent qualitative shifts in dynamics brought about by small changes in disturbance patterns, they can help predict tipping and regime changes in ecological and other systems. To study bifurcations in a disturbance parameter $\p\in\R$, we augment the domain of our disturbance rate function $r$.

\begin{definition}
    A parameterized disturbance rate function is a $C^2$ function $r(u,\p)$ from $U\times\R$ to $\R^n$ that is a disturbance rate function for each fixed $\p$.
\end{definition}
\noindent  
We study the connections between bifurcations in the ODE system
\[
x'=f(x)+r(x,\lambda)
\]
and those of the flow-kick map
\[
\Phi_\tau(x,\p)=\varphi_\tau(x)+\tau r(\vvarphi_\tau(\xx),\p).
\]
\noindent We will denote the first argument of $r$ as $x$ in the continuous context and as $u$ in the flow-kick context; the latter is for ease of tracking composition of $r$ with $\varphi$.  We are interested in whether a flow-kick model shares the bifurcations of its continuous analog, in which case the simpler ODE modeling framework might suffice to capture big changes in dynamics due to small changes in the disturbance rate parameter $\p$.
We first treat saddle-node bifurcations, then transcritical. In Sections \ref{sec:examples} we consider extensions to Hopf bifurcations in the context of a predator-prey example. 

\subsection*{Saddle-node bifurcations}
We restrict our attention to saddle-node bifurcations in the parameter $\p$ for the vector field $f(x)+r(x,\p)$ that satisfy the hypotheses of Sotomayor's Theorem  \cite{perko2001}. Table \ref{table:SNnew} provides these five sufficient conditions for a saddle-node bifurcation in the center column and the parallel conditions for the flow-kick map $\Phi_{\tau}(x,\p)$ in the right column \cite{robinson1999dynamical}. Our next theorem establishes that such saddle-node bifurcations continue locally from ODE to flow-kick systems, suggesting that ODEs suffice to capture these structures for small and frequent disturbances.

\begin{table}[b!]
\centering
\renewcommand{\arraystretch}{1.5}
\begin{tabularx}{6.1in}{XXX}
     & Vector field $f(x)+r(x,\p)$ \hspace*{0.5cm}& Map $\Phi_\tau(x,\p)$ \\
     \hline
     Saddle-node bifurcation at & $(x_*,\p_*)$ & $(x_\tau,\p_\tau)$\\
    \hline
     (i) Fixed point & $f(x_*)+r(x_*,\p_*)=0$ & $\Phi_\tau(x_\tau,\p_\tau)=x_\tau$\\
     (ii) Degenerate eigenvalue & \RaggedRight{$D_x[f+r](x_*,\p_*)$ has eigenvalue 0 with left, right eigenvectors $w_*,v_*$}& \RaggedRight{$D_x\Phi_\tau(x_\tau,\p_\tau)$ has eigenvalue 1 with left, right eigenvectors $w_\tau,v_\tau$}\\
     (iii) Other eigenvalues... & ...of \RaggedRight{$D_x[f+r](x_*,\p_*)$ have nonzero real parts} & \RaggedRight{... of $D_x\Phi_\tau(x_\tau,\p_\tau)$ have modulus $\not=1$}\\
     (iv) Transversality & $w_*\frac{\partial}{\partial \p}[f+r](x_*,\p_*) \not=0$ & $w_\tau\frac{\partial \Phi_\tau}{ \partial \p}(x_\tau,\p_\tau)\not=0$\\
     (v) Quadratic dominance & {\footnotesize $w_*(D_x^2[f+r](x_*,\p_*)(v_*,v_*))  \not=0$} & {\small $w_\tau (D_x^2\Phi_\tau(x_\tau,\p_\tau)(v_\tau,v_\tau))\not=0$}\\
     \hline
\end{tabularx}
\caption{Sufficient conditions for saddle-node bifurcations in $\p$ for vector fields \cite{perko2001} and maps \cite{robinson1999dynamical}.}
\label{table:SNnew}
\end{table}

\begin{theorem}\label{thm:s-n}
If the system $x'=f(x)+r(x,\p)$ ($x\in\R^n$) satisfies Sotomayor's conditions for a saddle-node bifurcation in the parameter $\p$ at $(x_*,\p_*)$ given in the center column of Table \ref{table:SNnew}, then for sufficiently small $\tau>0$ there exists a branch $(x_\tau,\p_\tau)$ of saddle-node bifurcations in $\p$ for the flow-kick maps $\Phi_\tau(x,\p)=\varphi_\tau(x)+\tau r(\vvarphi_\tau(\xx),\p)$ that varies continuously in $\tau$ and limits to $(x_*,\p_*)$ as $\tau\to 0$. 
\end{theorem}

\begin{proof}
It suffices to confirm the five conditions for a map to undergo a saddle-node bifurcation, given in the rightmost column of Table \ref{table:SNnew} \cite{robinson1999dynamical}. We first use the Implicit Function Theorem to construct a local branch of fixed points $(x_\tau,\lambda_\tau)$ with  eigenvalue 1, satisfying conditions (i) and (ii), then confirm that nondegeneracy conditions (iii), (iv), and (v) carry over from the continuous system locally along this branch.

Without loss of generality assume that the fixed point $x_*=0$; otherwise change coordinates by translating $x$ to $x-x_*$. Also without loss of generality, assume that the left and right eigenvectors of $D_x[f+r](x_*,\p_*)$ are $w_*=[1 \ 0 \cdots 0]$ and $v_*=[1 \ 0 \cdots 0]^T$, respectively; otherwise perform a linear change of basis that transforms $D_x[f+r](x_*,\p_*)$ into Jordan normal form with eigenvalue 0 in the first diagonal entry and matrix $B$ satisfying $\text{Re}(\text{spec} \ B)\not=0$ that sits minor to the $(1,1)$ position. 
Towards future economy in notation, denote the original system as
\begin{equation}\label{eq:orig_tilde}
    \tilde x'= \tilde f(\tilde x)+ \tilde r(\tilde x,\p).
\end{equation}
Let $f(\tilde x)=\tilde f(\tilde x)-D_{\tilde x}[\tilde f](\tilde x_*,\p_*)\tilde x$ and $r(\tilde x,\p)=\tilde r(\tilde x,\p)-D_{\tilde x}[\tilde r](\tilde x_*,\p_*)\tilde x$ represent the nonlinear parts of $\tilde f$ and $\tilde r$ at $(\tilde x_*,\p_*)$. Furthermore, 
let $x\in\R$ hold the first component of $\tilde x$ and $y\in\R^{n-1}$ hold the 2$^{\text{nd}}$ through n$^{\text{th}}$ components of $\tilde x$. We will adjust functions of $\tilde x$ in the natural way to be functions of $x$ and $y$ and reconstitute $\tilde x$ from its components using the notation $(x,y)=\tilde x$.
Equation \ref{eq:orig_tilde}
can then be expressed as
\begin{align}
    x'&=0 x+f_1(x,y)+ r_1(x,y,\p)\\
    y'&=By+f_{2{\text -}n}(x,y)+r_{2{\text -}n}(x,y,\p)
\end{align}
where subscripts indicate vector component(s) and $2\text{-}n$ means components $2$ through $n$. One can readily confirm that the following conditions hold at the fixed point $( x_*,y_*,\p_*)=(0,0,\p_*)$: 
\begin{align}
    f_1(x_*,y_*)+ r_1( x_*, y_*,\p_*)&=0,\label{eq:sn1}\\
    B y_* +  f_{2\text{-}n}( x_*, y_*)+r_{2\text{-}n}( x_*, y_*,\p_*)&=0,\label{eq:sn2}\\
   D_{(x,y)}[ f_1+  r_1]( x_*, y_*,\p_*)&=0, \label{eq:sn3}\\
   \text{and } D_{(x,y)}[ f_{2\text{-}n}+ r_{2\text{-}n}]( x_*, y_*,\p_*)&=0.\label{eq:sn4}
\end{align}

Let $\varphi(\tau,x,y)$ be the flow corresponding to the undisturbed system
\begin{align}
    x'&=0x+f_1(x,y)\\
    y'&=By+f_{2\text{-}n}(x,y),
\end{align}
so that $\Phi_\tau(x,y,\p)=\varphi_\tau(x,y)+\tau r(\phi_\tau(x,y),\p)$. We define a function $\hat F$ whose zeros give fixed points of the flow-kick map $\Phi_\tau$ with eigenvalue 1: 
\begin{align}
    \widehat F(\tau,x,y,\p)&=
    \begin{bmatrix}
    \varphi_1(\tau,x,y)+\tau r_1( \varphi_\tau(x,y),\p)-x\\
         \varphi_{2\text{-}n}(\tau,x,y)+\tau r_{2\text{-}n}( \varphi_\tau(x,y),\p)-y\\
    \det(D_{(x,y)}\Phi_\tau-1I)
    \end{bmatrix}.
\end{align}
Note that for $\tau=0$, $D_{(x,y)}\Phi_0=I$ and all choices of $x$, $y$, and $\p$ yield zeros of $\widehat F$. To eliminate this degeneracy, we will remove factor(s) of $\tau$ from each component of $\widehat F$. We use a Taylor expansion of the flow in $\tau$ as in equations \eqref{eq:taylormessy}-\eqref{eq:taylorclean} to get
\begin{align}
    \widehat F(\tau,x,y,\p)
    &=\begin{bmatrix}
        \tau (f_1(x,y)+r_1(\varphi_\tau(x,y),\p)+\mathcal{O}(\tau))\\
        \tau (By+f_{2\text{-}n}(x,y)+ r_{2\text{-}n}(\varphi_\tau(x,y),\p)+\mathcal{O}(\tau))\\
    \det(D_{(x,y)}\Phi_\tau-1I)
    \end{bmatrix} 
\end{align}
where $\mathcal{O}(\tau^k)$ denotes a term with a factor of $\tau^k$. We unpack the final component as follows:
\begin{align}
D_{(x,y)}\Phi_\tau
    &=D_{(x,y)}
    \begin{bmatrix}
        x+\tau f_1(x,y) +\mathcal{O}(\tau^2)+\tau r_1(\varphi_\tau(x,y),\p)\\
        y+\tau By +\tau f_{2\text{-}n}(x,y)+\mathcal{O}(\tau^2)+\tau r_{2\text{-}n}(\varphi_\tau(x,y),\p)
    \end{bmatrix}\label{eq:DxyPhitau}\\
    & = D_{(x,y)} \begin{bmatrix}x \\ y \end{bmatrix} + \tau D_{(x,y)}\begin{bmatrix} f_1(x,y)+\mathcal{O}(\tau)+r_1(\varphi_\tau(x,y),\p) \\
By + f_{2\text{-}n}(x,y)+\mathcal{O}(\tau) + r_{2\text{-}n}(\varphi_\tau(x,y),\p)\end{bmatrix} \\
    & = I + \tau D_{(x,y)} M
\end{align}
where 
\begin{equation}\label{eq:M}
M = \begin{bmatrix} f_1(x,y)+\mathcal{O}(\tau)+r_1(\varphi_\tau(x,y),\p) \\
By + f_{2\text{-}n}(x,y)+\mathcal{O}(\tau) + r_{2\text{-}n}(\varphi_\tau(x,y),\p)\end{bmatrix}.
\end{equation} 
Then 
\begin{align}
    \det(D_{(x,y)}\Phi_\tau - I) = \det (I+\tau D_{(x,y)}M - I) = \tau^n \det(D_{(x,y)} M)\label{eq:detDM}
\end{align}
and we choose $F_3(\tau,x,y,\p)=\det(D_{(x,y)}M)$ to eliminate a factor of $\tau^n$ from \eqref{eq:detDM}. Let
\begin{align}
    F(\tau,x,y,\p)&=\begin{bmatrix}
        f_1(x,y)+r_1(\varphi_\tau(x,y),\p)+\mathcal{O}(\tau)\\
        By+f_{2\text{-}n}(x,y)+ r_{2\text{-}n}(\varphi_\tau(x,y),\p)+\mathcal{O}(\tau)\\
    F_3(\tau,x,y,\p)
    \end{bmatrix} . 
\end{align}
Since zeros of $F$ are also zeros of $\widehat F$, they correspond to flow-kick fixed points with eigenvalue 1. We apply the Implicit Function Theorem to $F$. 

The fact that $F(0,x_*,y_*,\p_*)=0$ can be  confirmed from properties \eqref{eq:sn1}, \eqref{eq:sn2}, and \eqref{eq:sn3} of the continuous system (see Appendix \ref{sec:bif_app}, equations \eqref{eq:F3_start}-\eqref{eq:F3_end} for computation of $F_3$).

Next we check the invertibility of $F$ at $(\tau,x,y,\p)=(0,x_*,y_*,\p_*)$. At $(\tau,x,y,\p)=(0,x_*,y_*,\p_*)$ the  partial derivatives of $f$ and $r$ with respect to $x$ and $y_j$ vanish and $D_{(x,y,\p)}F$ evaluates to
\begin{equation}\label{eq:block}
\left[
\begin{array}{c|ccc|c}
    0 & \cdots & 0 & \cdots & \frac{\partial r_1}{\partial \p}(x_*,y_*,\p_*)\\
\hline
    \vdots &  &  & & \vdots \\
    0 &  & B & & *\\
    \vdots & & & & \vdots\\
\hline
    \frac{\partial F_3}{\partial x}(0,x_*,y_*,\p_*) & \cdots & * & \cdots & *
\end{array}   
\right].
\end{equation}
Because $B$ is invertible, the matrix \eqref{eq:block} is invertible if and only if both $\frac{\partial r_1}{\partial \p}(x_*,y_*,\p_*)\not=0$ and $\frac{\partial F_3}{\partial x}(0,x_*,y_*,\p_*)\not=0$. 
We confirm both conditions in Appendix \ref{sec:bif_sn}, using, respectively, the transversality and quadratic dominance conditions in the continuous system. 

One can check that $F$ is $C^1$, so by the Implicit Function Theorem,
for $\tau$ sufficiently small, there exists a unique branch $(x_\tau, y_\tau, \p_\tau)$ of solutions to $F(\tau, x,y,\p)=0$ that varies continuously in $\tau$ from $(x_0,y_0,\p_0)=(x_*,y_*,\p_*)$. By construction, $(x_\tau,y_\tau)$ are fixed points with eigenvalue $1$ for the flow-kick maps $\Phi_\tau(x,y,\p)=\varphi_\tau(x,y)+\tau r(\varphi_\tau(x,y),\p_\tau)$.

The nondegeneracy conditions (iii), (iv), and (v) of the flow-kick map follow from continuity arguments. 
We identify expressions within each of the map conditions (iii), (iv), and (v) that are continuous functions of $\tau$ and are nonzero when $\tau = 0$. From these we conclude that each condition holds locally, so there exists an interval $\tau\in (0,\delta)$ such that all conditions are satisfied simultaneously. Details are provided in Appendix \ref{sec:bif_sn}.
\end{proof}

Theorem \ref{thm:s-n} concerns bifurcations in the parameter $\p$. However, we also saw a saddle-node bifurcation in $\tau$ for the logistic system in Figure \ref{fig:logbifnew}c and Figure \ref{fig:cor1}, with $\p=-0.24$ fixed. 
One can interpret this bifurcation in $\tau$ as a sign that if the average harvest rate is held constant, discrete harvests that were sustainable at high frequencies and low amplitudes become unsustainable at lower frequencies and higher amplitudes.

\subsection*{Transcritical bifurcations}

Since transcritical bifurcations are not generic, we expect that additional conditions are needed to ensure a transcritical bifurcation continues from ODE to flow-kick systems.

A natural condition to consider is that the disturbance $r(x,\p)$ is zero at a fixed point $x_*$ for all values of the parameter $\p$. This occurs, for example, in the multiplicative disturbance setting where $x_*=0$ and $r(x,\p)=\p x$. The following theorem combines this hypothesis with Sotomayor's characterization of transcritical bifurcations to demonstrate conditions under which these bifurcations continue from ODE to flow-kick models.

\begin{theorem}\label{thm:tc}
    Suppose that at $(x_*,\p_*)$ the ODE $x'=f(x)+r(x,\p)$ satisfies Sotomayor's conditions for a transcritical bifurcation \cite{perko2001}, given in the middle column of Table \ref{table:TC}. Suppose further that $r(x_*,\p)=0$ for all $\p$. Then a branch $(x_*,\p_\tau)$ of transcritical bifurcations in $\p$ continues locally for the flow-kick maps $\Phi_\tau(x,\p)=\varphi_\tau(x)+\tau r(\varphi_\tau(x),\p)$, varying continuously in $\tau$ and limiting to $(x_*,\p_*)$ as $\tau\to 0$.
\end{theorem}

\begin{table}[b!]
\centering
\renewcommand{\arraystretch}{1.5}
\begin{tabularx}{6.1in}{XXX}
    & Vector field $f(x)+r(x,\p)$ \hspace*{0.5cm}& Map $\Phi_\tau(x,\p)$ \\
     \hline
     Transcritical Bif'n. at & $(x_*,\p_*)$ & $(x_*,\p_\tau)$\\
    \hline
     (i) Fixed point & $f(x_*)+r(x_*,\p_\tau)=0$ & $\Phi_\tau(x_*,\p_\tau)=x_*$\\
     (ii) Degenerate eigenvalue & \RaggedRight{$D_x[f+r](x_*,\p_*)$ has simple eigenvalue 0 with left, right eigenvectors $w_*,v_*$}& \RaggedRight{$D_x\Phi_\tau(x_*,\p_\tau)$ has simple eigenvalue 1 with left, right eigenvectors $w_\tau,v_\tau$}\\
     (iii) Other eigenvalues... & ...of \RaggedRight{$D_x[f+r](x_*,\p_*)$ have nonzero real parts} & \RaggedRight{... of $D_x\Phi_\tau(x_*,\p_\tau)$ have modulus $\not=1$}\\
     (iv) Non-transversality & $w_*\frac{\partial}{\partial \p}[f+r](x_*,\p_*) =0$ & $w_\tau\frac{\partial \Phi_\tau}{ \partial \p}(x_*,\p_\tau)=0$\\
     (v) Mixed partials & {\footnotesize $w_*\left(D_x\frac{\partial}{\partial \p}[f+r](x_*,\p_*)v_*\right)\not=0$} & $w_\tau\left(D_x\frac{\partial \Phi_\tau}{ \partial \p}(x_*,\p_\tau)v_\tau\right)\not=0$\\
     (vi) Quadratic dominance & {\footnotesize $w_*(D_x^2[f+r](x_*,\p_*)(v_*,v_*))  \not=0$} & {\small $w_\tau (D_x^2\Phi_\tau(x_*,\p_\tau)(v_\tau,v_\tau))\not=0$}\\
     \hline
\end{tabularx}
\caption{Sufficient conditions for transcritical bifurcations for vector fields \cite{perko2001}, and analogous conditions for maps.}
\label{table:TC}
\end{table}

\begin{proof}
Note that since $ f( x_*)+ r( x_*,\p_*)=0$ and $r( x_*,\p)=0$ for all $\p$, we have $f( x_*)=0$ and $\varphi_\tau( x_*)= x_*$. It follows that $ x_*$ is a fixed point for the flow-kick map $\Phi_\tau( x,\p)=\varphi_\tau( x)+\tau r( \varphi_\tau(x),\p)$ for all values of $\p$. 

Next we construct a local branch of parameter values $\p_\tau$ that yield an eigenvalue of $1$ at the fixed point $x_*$ for each flow-kick map $\Phi_\tau$. As in the proof of Theorem \ref{thm:s-n}, change coordinates to rewrite the original ODE $\tilde x'=\tilde f(\tilde x)+\tilde r(\tilde x,\p)$ as 
\begin{align*}
    x'&=0 x+f_1(x,y)+ r_1(x,y,\p)\\
    y'&=By+f_{2{\text -}n}(x,y)+r_{2{\text -}n}(x,y,\p)
\end{align*}
where $x\in\R$, $y\in\R^{n-1}$, $B$ is hyperbolic, and spatial derivatives of $f_i+r_i$ are zero at $(x_*,y_*,\p_*)$ (see equations \eqref{eq:orig_tilde}-\eqref{eq:sn4}). 
We seek roots of 
\begin{equation}\label{eq:detDxPhi-I}
    \det\left(D_{\tilde x}\Phi_\tau(\tilde x_*,\p)-1I\right)
    =\tau^n \det(D_{(x,y)}M)\big|_{(x_*,y_*)}
\end{equation} 
where 
\begin{equation}
M = \begin{bmatrix} f_1(x,y)+\mathcal{O}(\tau)+r_1(\varphi_\tau(x,y),\p) \\
By + f_{2\text{-}n}(x,y)+\mathcal{O}(\tau) + r_{2\text{-}n}(\varphi_\tau(x,y),\p)\end{bmatrix}
\end{equation}
(see equations \eqref{eq:DxyPhitau}-\eqref{eq:M}). Setting up for the Implicit Function Theorem, we eliminate the factor $\tau^n$ from \eqref{eq:detDxPhi-I} and define 
{\allowdisplaybreaks
\begin{align}
    \begin{split}
    F(\tau,\p)&=\det(D_{(x,y)}M)\big|_{(x_*,y_*)}\\
    &=\Bigg[ \left(\frac{\partial f_1}{\partial x} + \frac{\partial r_1}{\partial x} +\frac{\partial}{\partial x}\mathcal{O}(\tau) \right)\det\left(B+D_y\left[f_{2\text{-}n}+r_{2\text{-}n} \right] + D_y\mathcal{O}(\tau) \right) 
\end{split}\\
\begin{split}
    & \hspace{1cm} +\sum\limits_{j=1}^{n-1}(-1)^j\left(\frac{\partial f_1}{\partial y_j}+ \frac{\partial r_1}{\partial y_j} +\frac{\partial}{\partial y_j}\mathcal{O}(\tau) \right) \det(M_{1,j+1} )\Bigg]\Bigg|_{(x_*,y_*)}
\end{split}\\
&=\alpha(\tau,x_*,y_*,\p) \, \beta(\tau,x_*,y_*,\p)+\sum\limits_{j=1}^{n-1}(-1)^j\gamma_j(\tau,x_*,y_*,\p)\det(M_{1,j+1})
\end{align}
}
(see equations \eqref{eq:rearr}-\eqref{eq:det_end} for the derivative and determinant computations and \eqref{eq:alpha}-\eqref{eq:gamma} for explicit definitions of $\alpha$, $\beta$, and $\gamma$).

Zeros of $F$ are also zeros of $\det\left(D_{\tilde x}\Phi_\tau(\tilde x_*,\p)-1I\right)$ and represent $(\tau,\p)$ parameter combinations that yield an eigenvalue of 1 for $D_{\tilde x}\Phi_\tau$ at $\tilde x_*=(x_*,y_*)$. We have that $F(0,\p_*)=0$ because spatial derivatives of $f_1+r_1$ are zero at $(x_*,y_*,\p_*)$ (see equation \eqref{eq:sn3}). To extend $\p$ as a function of $\tau$ satisfying $F(\tau,\p_\tau)=0$, we require that $\partial F/\partial \p(0,\p_*)\not=0$. We have
\begin{align}
    \frac{\partial F}{\partial \p}
    &=\frac{\partial \alpha}{\partial \p} \, \beta +\alpha\frac{\partial \beta}{\partial \p} +\sum\limits_{j=1}^{n-1}(-1)^j\frac{\partial \gamma}{\partial \p}\det(M_{1,j+1}) 
    +\sum\limits_{j=1}^{n-1}(-1)^j\,\gamma\,\frac{\partial}{\partial \p}\big[\det(M_{1,j+1})\big].
\end{align}
As noted in Appendix \ref{sec:bif_sn}, at $(\tau,x,y,\p)=(0,x_*,y_*,\p_*)$ the terms $\alpha$ and $\gamma$ containing spatial derivatives of $f_1+r_1$ vanish, while $\beta$ reduces to $\det B$. Further, $\det M_{1,j+1}=0$. Thus   
\begin{align}    
    \frac{\partial F}{\partial \p}(0,\p_*)
    &=\det(B) \ \frac{\partial \alpha}{\partial \p}\Bigg|_{(0,x_*,y_*,\p_*)} \\
    &=\det(B) \ \frac{\partial}{\partial \p}\left[ \frac{\partial f_1}{\partial x} + \nabla_u r_1\!\cdot\!\frac{\partial \varphi_\tau}{\partial x} +\frac{\partial}{\partial x}\mathcal{O}(\tau)\right]\Bigg|_{(0,x_*,y_*,\p_*)}\\
    &=\det(B) \ \left(\frac{\partial \, \nabla_u r_1}{\partial \p}\!\cdot\!\frac{\partial \varphi_0}{\partial x} \right)\Bigg|_{(x_*,y_*,\p_*)}\\
    &=\det(B) \ \left[\frac{\partial^2 r_1}{\partial\p\partial x}(x_*,y_*,\p_*) \ \cdots \ \frac{\partial^2 r_1}{\partial\p\partial y_{n-1}}(x_*,y_*,\p_*)\right] \ \left[1 \ 0 \, \cdots \, 0 \right]^T\\
    &=\det(B) \ \frac{\partial^2 r_1}{\partial\p\partial x}(x_*,y_*,\p_*). 
\end{align}
It follows from the mixed partials condition (v) on the vector field $\tilde f + \tilde r$ (see Table \ref{table:TC}) and Clairaut's Theorem that $\frac{\partial^2 r_1}{\partial\p \partial x}(x_*,y_*,\p_*)\not=0$ (see Appendix \ref{sec:bif_tc}). And from the assumption $\text{Re}(\text{spec} \ B)\not=0$ we get $\det(B)\not=0$. Hence $\frac{\partial F}{\partial \p}(0,\p_*)\not=0$.
Lastly, one can confirm that $F(\tau,\p)$ is $C^1$. The Implicit Function Theorem then implies there exists $\p_\tau$,  a local continuous function of $\tau$, such that $\p_0=\p_*$ and $F(\tau,\p_\tau)=0$. By construction, the points $(x_*,\p_\tau)$ are fixed points for the flow-kick maps $\Phi_\tau$ and $D_{\tilde x}\Phi_\tau(x_*,\p_\tau)$ has eigenvalue 1.

To confirm the non-transversality condition (iv), we compute in original coordinates ($x\in\R^n$)
\begin{align}
    \frac{\partial \Phi_\tau}{\partial \p}(x_*,\p_\tau)&=\frac{\partial }{\partial \p}\Big[\varphi_\tau(x)+\tau r(\varphi_\tau(x),\p) \Big]\Big|_{(x_*,\p_\tau)}\label{eq:nt1}\\
    &=\tau\frac{\partial r}{\partial \p}(\varphi_\tau(x_*),\p_\tau)\label{eq:nt2}\\
    & = \tau \frac{\partial r}{\partial \p}(x_*,\p_\tau)\\
    &=0
\end{align}
where the final equality follows from the fact that $r(x_*,\p)=0$ for all $\p$.

The conditions for non-degenerate eigenvalues (iii) and for quadratic dominance (vi) are identical to those in Theorem \ref{thm:s-n}. Additionally, the mixed partials condition (v) is nonzero because of a similar continuity argument. The details are provided in Appendix \ref{sec:bif_tc}. 
\end{proof}

\section{Examples}\label{sec:examples}

In this section, we use models of vegetation and precipitation dynamics (Subsection \ref{sec:klaus}) and predator-prey interactions (Subsection \ref{sec:predprey}) to illustrate how fixed points and their stabilities, as well as bifurcations, continue from ODE to flow-kick systems. Additionally, we numerically compute continuations of a Neimark-Sacker bifurcation, offering evidence that Hopf bifurcations in continuous systems may continue to Neimark-Sacker bifurcations more generally.

Numerical continuations of fixed points and bifurcations are calculated using MatContM \cite{dhooge2008new}. To compute fixed point and bifurcation continuations numerically, MatContM requires first, second, and third derivatives of the map \cite{meijer2017matcontm}. Since we do not have a map in closed form, we bypass the autodifferentation algorithm in MatContM and encode these derivatives using central differences on perturbed initial conditions. 
This is a higher error method than autodifferentation, but for these two-dimensional systems, the continuations are successful.

\subsection{Nonspatial Klausmeier model of vegetation and precipitation}\label{sec:klaus}

In \cite{klausmeier1999regular}, Klausmeier briefly considers a nonspatial model of vegetation and water in semi-arid conditions. Rainfall appears as a small continuous disturbance to the water variable,
representing average precipitation rates (total precipitation divided by time). The model can be regarded as phenomenological, relating average biomass to rainfall input on the timescale of years. Of course, actual precipitation events occur on much faster timescales than vegetation growth. Here we illustrate how fixed points and bifurcations change when we instead model rain as a periodic kick. This example serves more for mathematical illustration---and less for biological prediction---as the flow phase of the flow-kick model may not capture vegetation dynamics under pulsed rainfall.

\begin{figure}
    \centering
    \includegraphics[width=.9\textwidth]{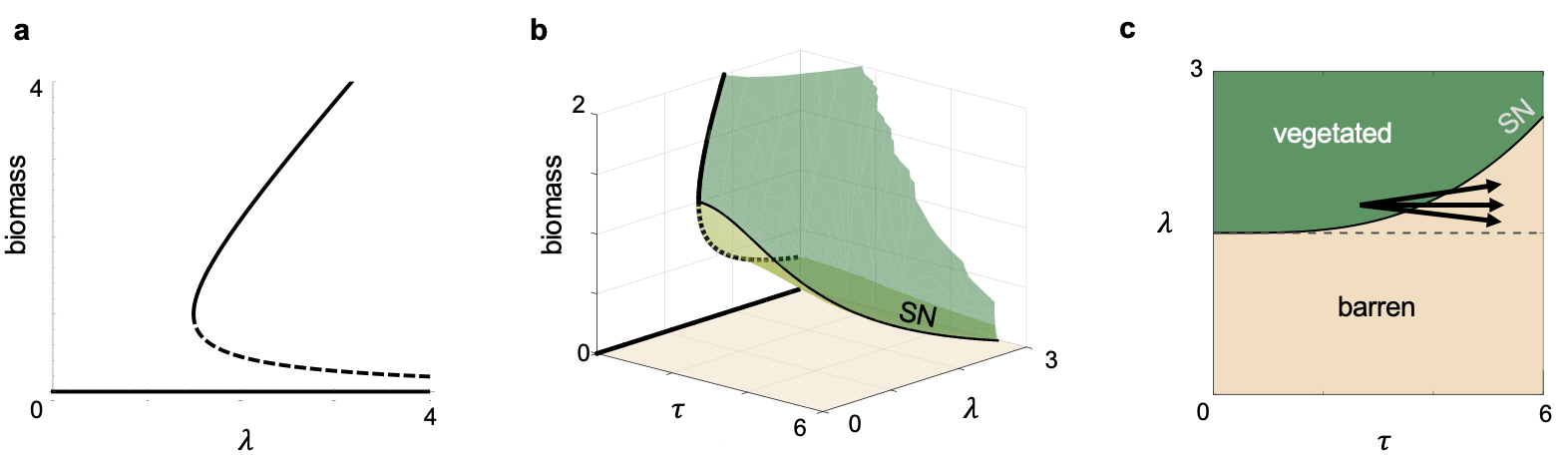}
    \caption{\small{(a) Equilibrium biomass in the continuous nonspatial Klausmeier model \eqref{eq:klaus} with $m=0.75$, which yields a saddle-node bifurcation at $\p=1.5$. Solid and dashed lines represent stable and unstable branches of equilibria, respectively. (b) Continuation of fixed points and saddle-node bifurcation from continuous rain model \eqref{eq:klaus} (in $\tau=0$ plane) to discrete rain model \eqref{eq:klaus_flow-kick} (with $\tau>0$ and $\p$ variable). Dark green and tan surfaces contain stable fixed points; light green contains unstable fixed points. (c) By projecting  (b) down the biomass axis onto the $(\tau,\p)$ disturbance parameter space, we obtain a stability diagram.}}
\label{fig:klaus_continuation}
\end{figure}

In the nonspatial Klausmeier model, let $x=(x_1,x_2)$ where $x_1$ represents nondimensional vegetation biomass and $x_2$ represents nondimensional water. The parameter $m$ controls vegetation mortality rates, and $\p$ 
represents water input via precipitation:
\begin{equation}
    \label{eq:klaus}
    \begin{aligned}
       x_1' & = x_2x_1^2-mx_1\\
       x_2' & = -x_2x_1^2-y +\p.
    \end{aligned}
\end{equation}
Klausmeier shows that this system has one stable vegetation-free (barren) equilibrium at $(0,\p)$ for all values of $\p$ and $m$. When $\p=2m$, the system undergoes a saddle node bifurcation at $(1,\p/2)$, so that for $\p \geq 2m$ the system has two vegetated equilibria---one stable, the other unstable---at $(\frac{\p\pm \sqrt{\p^2-4m^2}}{2m}, \frac{\p\mp \sqrt{\p^2-4m^2}}{2})$. 
Figure \ref{fig:klaus_continuation}a shows stable (solid) and unstable (dashed) branches of  equilibrium biomass for $m=0.75$ and $\p$ varying.

In contrast, we consider a flow-kick model of precipitation with the flow generated by the rain-free ODE
\begin{equation}
    \label{eq:klaus_undisturbed}
    \begin{aligned}
        x_1' & = x_2x_1^2-mx_1\\
        x_2' & = -x_2x_1^2-x_2.
    \end{aligned}
\end{equation}
The flow-kick map is then 
\begin{equation}
\label{eq:klaus_flow-kick}
    \Phi_{\tau}(x,\p) = \varphi_{\tau}(x)+ \tau   \begin{bmatrix}
    0 \\ \p
\end{bmatrix}
\end{equation} where $\varphi_{\tau}(x)$ solves \eqref{eq:klaus_undisturbed}.

By Theorem \ref{thm:stability_n},  hyperbolic equilibria of the continuous rainfall model \eqref{eq:klaus} continue locally to fixed points of the discrete rainfall model \eqref{eq:klaus_flow-kick}. Figure \ref{fig:klaus_continuation}b illustrates  biomass equilibrium branches for \eqref{eq:klaus} in the $\tau=0$ plane as well as numerical continuations of vegetated flow-kick fixed points for $\tau>0$. A vertical section of this curved surface with fixed $\p$ traces out branches of flow-kick fixed points that continue from the ODE equilibria, except at the bifurcation point $\p=1.5$. In addition, one can analytically confirm that each vegetation-free equilibrium for \eqref{eq:klaus} at $(0,\p)$ continues for $\tau>0$ to vegetation-free flow-kick fixed points at $(0,\tau\p/(1-e^{-\tau}))$. 

Although flow-kick fixed points do not continue from the nonhyperbolic vegetated equilibrium with $\p=1.5$ fixed, the saddle-node bifurcation that occurs here in the ODE system \eqref{eq:klaus} continues numerically in $\tau$ with $\p$ variable to the curves labeled SN in Figure \ref{fig:klaus_continuation}bc. 

It is also of note that the precipitation rate $\p$ at which the saddle-node bifurcation occurs grows larger as the period between precipitation events $\tau$ increases. In other words, higher precipitation averages are needed to maintain vegetation when rain arrives in more intense, less frequent bursts. 
Within the model, increases in the duration of droughts---with or without changes to average precipitation rates---could tip a system over the bifurcation curve from a vegetated to barren regime, as illustrated by the black arrows in Figure \ref{fig:klaus_continuation}c. Caution is warranted in drawing biological predictions from such a simplified model, but these results do align with predictions of more sophisticated models, discussed next.

Precipitation events can be treated more realistically in a stochastic framework---for example, as a Poisson point process with exponentially distributed depths \cite{rodriguez1999probabilistic}. Gandhi and colleagues \cite{gandhi2023pulsed} recently coupled this approach with a spatially explicit reaction-diffusion model of vegetation and water dynamics. They found that the choice of deterministic versus stochastic rainfall events altered qualitative outcomes in some parameter regimes of the model, including the appearance of vegetation bands and their direction of migration. 
On the other hand, some insights from their stochastic framework align with our deterministic predictions. In particular, in scenarios of identically low mean annual precipitation rates, vegetation collapsed more quickly on average as the mean storm depth increased and mean arrival rates of storms decreased. 
Thus despite its limitations, a deterministic and non-spatial model of rainfall that resolves rain events in time can reveal patterns consistent with more realistic models.

The spatial and stochastic aspects of disturbance modeling that we have ignored for much of this paper are interesting directions for future research. For example, is it possible to predict general effects of discretizing disturbances, starting from a partial differential equation or stochastic differential equation model?

\subsection{Predator-Prey Model}\label{sec:predprey}  

We use a predator-prey model to illustrate additional future directions for this work, including the continuation of periodic orbits and Hopf bifurcations.

The following ODE system models undisturbed prey ($x$) and predator ($y$) populations \cite{may1972limit,rosenzweig1971paradox}, with parameters chosen to yield a stable limit cycle in the first quadrant:

\begin{equation}\label{eq:predprey_undisturbed}
        \begin{aligned}
            \frac{dx}{dt}&=x(1-x/4)-0.5y(1-e^{-1.5x})\\
            \frac{dy}{dt}&=-0.5y+y(1-e^{-0.5x}).
        \end{aligned}
    \end{equation}
Consider disturbance to \eqref{eq:predprey_undisturbed} in the form of proportional harvesting of the predator ($y$). On the one hand, a continuous model of this disturbance is
\begin{equation}\label{eq:predprey}
        \begin{aligned}
            \frac{dx}{dt}&=x(1-x/4)-0.5y(1-e^{-1.5x})\\
            \frac{dy}{dt}&=-0.5y+y(1-e^{-0.5x})-\p y,
        \end{aligned}
\end{equation}
where $\p$ controls the harvesting rate. On the other hand, a discrete model of this disturbance is given by the flow-kick map
\begin{equation}\label{eq:predprey_flowkick}
    \Phi_{\tau }(x,y,\lambda) = \varphi_{\tau}(x,y)- \tau\begin{bmatrix}
        0 \\ \p [\varphi_{\tau}(x,y)]_2
    \end{bmatrix}
\end{equation}
where $\varphi_{\tau}(x,y)$ solves \eqref{eq:predprey_undisturbed} and $[ \ \cdot \ ]_2$ takes its second component.

\begin{figure}
    \centering
\includegraphics[width=0.95\textwidth]{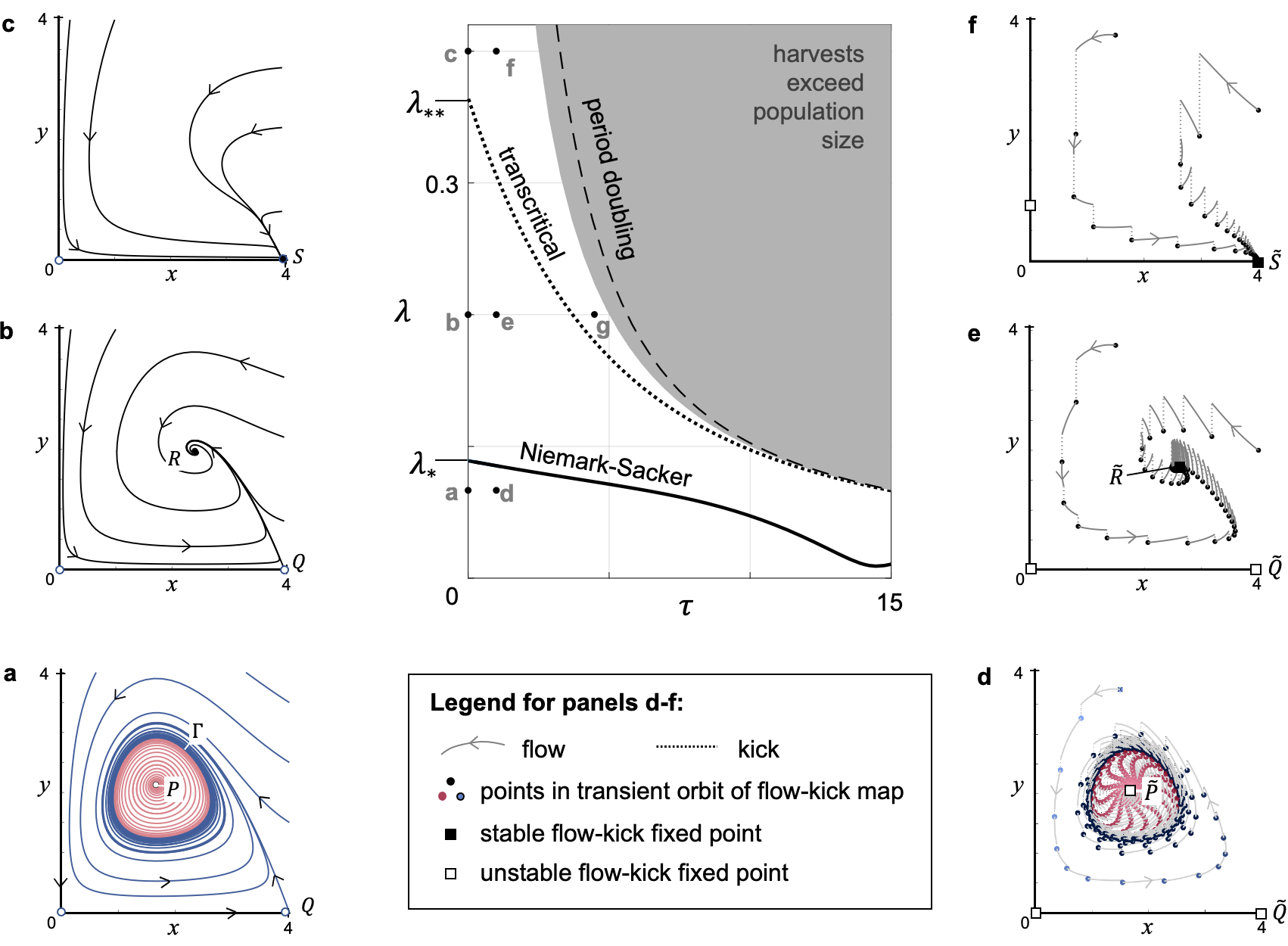}
    \caption{\small Predator prey dynamics under continuous and discrete disturbance harvests of predators. (a)-(c) Phase diagrams for the continuous model (visualized at $\tau = 0$ in the central panel), showing three qualitative states of the system for disturbance rates $\p=1/15$ (panel a), $\p=1/5$ (panel b), and $\p=2/5$ (panel c). (d)-(f) Phase diagrams for flow-kick systems at disturbance rates matching (a)-(c), respectively, with $\tau=1$. Qualitative behaviors are similar between continuous and discrete disturbance models. Center panel: 
    Continuation of transcritical and Hopf bifurcations in continuous disturbance model to transcritical  (dotted line) and Niemark-Sacker (solid line) bifurcations in the flow-kick disturbance model. The flow-kick system also exhibits period doubling bifurcations (dashed line) for nonphysical parameter combinations that cause harvest size to exceed population size (grey shaded area). 
    Note that this diagram may be incomplete as a stability diagram, because the flow-kick system may have additional bifurcations that are disconnected from the $\tau=0$ axis.
    } \label{fig:predprey_stability}
\end{figure}

We compare dynamics and bifurcations across the continuous and discrete disturbance models. Harvesting pressure drives three qualitative regimes in the continuous system \eqref{eq:predprey}. As illustrated in Figure \ref{fig:predprey_stability}a for $\p=1/15$, low values of the harvesting parameter preserve a stable limit cycle $\Gamma$ that encloses an unstable coexistence equilibrium $P$. As $\p$ increases, the limit cycle contracts and undergoes a supercritical Hopf bifurcation at $\p_* \approx 0.089$, yielding a stable coexistence equilibrium $R$ in the first quadrant (Figure \ref{fig:predprey_stability}b). As $\p$ increases further, this equilibrium moves towards the saddle at $Q=(4,0)$. The two equilibria collide in a transcritical bifurcation at $\p_{**}\approx 0.365$, and for larger harvesting rates the point $Q$ is stable (Figure \ref{fig:predprey_stability}c). Thus the continuous model predicts either cycles of predator and prey populations, steady coexistence, or crash of the predators, depending on harvesting rates.

The discrete disturbance model yields similar dynamics and predictions for small recovery times $\tau$. The rightmost panels d, e, and f in Figure \ref{fig:predprey_stability} give flow-kick phase portraits for $\tau=1$ and values of $\p$ that match the continuous systems depicted in the leftmost panels a, b, and c, respectively. As Theorem 1 predicts, the hyperbolic equilibria at $P$, $Q$, $R$, and $O$ (the origin) in the continuous harvesting model continue to flow-kick fixed points with the same stabilities at $\tilde P$, $\tilde Q$, $\tilde R$, and $\tilde O$ in the discrete harvesting model. One subtle difference in dynamics is best understood by returning to an impulsive ODE framework: when we account for distinct flow and kick phases, a flow-kick fixed point does not necessarily represent unchanging populations. For example, the fixed point $\tilde R$ occurs where transient predator growth balances harvests, and both nontrivial processes occur during one iteration of the flow-kick map. On the other hand, because the origin and $Q=\tilde Q=(4,0)$ are invariant under the flows generated by \eqref{eq:predprey_undisturbed} and \eqref{eq:predprey} as well as the kick to $y$, these points are both flow-kick fixed points and true equilibria of the impulsive DE across all values of $\tau$ and $\lambda$.

The predator-prey example also illustrates continuation of transcritical bifurcations. The central panel in Figure \ref{fig:predprey_stability} plots bifurcation curves for the flow-kick map in $(\tau,\p)$ parameter space. We developed the condition
\begin{equation}
    \p\tau=1-\text{exp}(\tau(\tfrac12-e^{-2}))
\end{equation}
necessary for a transcritical bifurcation to occur at $(x,y)=(4,0)$ (see Appendix \ref{sec:predprey_bifconditions}), which we numerically confirmed in MatContM with error less than $10^{-6}$. 
The dashed curve represents this condition and connects to the ODE transcritical bifurcation at $(0,\p_{**})$.

In addition to the similarities predicted by Theorem \ref{thm:stability_n} and \ref{thm:tc}, the numerical simulations of predator-prey dynamics in Figure \ref{fig:predprey_stability} display several patterns that warrant further study. 
First, the fixed points that continue retain not only their stabilities but also their type---for example, node, spiral, or saddle. We conjecture that equilibrium type continues in general from ODE to flow-kick systems when the type in the ODE system is not degenerate. 
Second, the stable limit cycle $\Gamma$ in Figure \ref{fig:predprey_stability}a has a striking counterpart in Figure \ref{fig:predprey_stability}d:  flow-kick trajectories tend towards what appears to be an invariant closed curve near $\Gamma$. Our numerical simulation does not distinguish high-period periodic points from a dense orbit along an invariant curve, but in either case the cycling of predator and prey species in the flow-kick system closely mimics that of the continuous system. Continuation of periodic orbits from ODEs to analogous structures in flow-kick models represents another area for future study. Third, we used MatcontM to identify a curve of Niemark-Sacker bifurcations (solid black line), which connects to the continuous Hopf bifurcation at $(0,\p_*)$. Given that Hopf bifurcations are generic, we conjecture that few additional conditions will be needed to guarantee that they continue from ODE systems to Niemark-Sacker bifurcations in flow-kick systems.

Despite the continuation patterns that the predator-prey system illustrates for small $\tau$, larger $\tau$ values can destroy qualitative similarities between ODE and flow-kick disturbance models. For example, even though $\p=0.2$ predicts stable coexistence between predators and prey when $\tau=1$, increasing the time between harvests to $\tau=4.5$ changes the prediction to predator crash (point $g$ in Figure \ref{fig:predprey_stability}, central panel). In addition, period doubling bifurcations occur in the flow-kick system at the point $(x,y)=(4,0)$ (determined numerically in MatContM and analytically, see Appendix \ref{sec:predprey_bifconditions}). Although they occur for ecologically unfeasible parameter values in which greater than 100\% of the predator population is harvested, they highlight the possibility that discretizing disturbance can introduce bifurcation structures not present in the analogous ODE.

\section{Discussion}\label{sec:disc}

We have presented both rigorous arguments and numerical evidence that ODE and flow-kick disturbance models resemble one another when the kicks occur at high frequencies. On the one hand, flow-kick maps $\Phi_\tau$ with the same average disturbance rate generate a vector field with that disturbance rate in the limit as $\tau$---the period between kicks---decreases to zero (Proposition \ref{prop:ctslimit}). On the other hand, we have seen multiple dynamic features of ODE disturbance models continue locally to flow-kick disturbance models. These include fixed points with their stabilities (Theorem \ref{thm:stability_n} and Figures \ref{fig:klaus_continuation} and \ref{fig:predprey_stability}), saddle node bifurcations (Theorem \ref{thm:s-n} and Figure \ref{fig:klaus_continuation}), and transcritical bifurcations (Theorem \ref{thm:tc} and Figure \ref{fig:predprey_stability}). Future work in this direction could treat continuation of periodic orbits and Hopf bifurcations in ODE disturbance models to analogous structures in flow-kick systems, as observed in Figure \ref{fig:predprey_stability}. These proven and conjectured patterns point towards a consistent take-away: if true disturbances are discrete and high-frequency, but we choose to model them with an ODE instead, then dynamic structures in the ODE model approximate similar structures in the kicked system. 

However, this story changes as disturbance frequency decreases. Although saddle-node and transcritical bifurcations do persist from continuous to flow-kick systems, 
they can occur at different disturbance rates. 
This phenomenon was previously observed for additively harvested populations with logistic growth \cite{zeeman2018resilience} and with an Allee effect \cite{meyer2018quantifying}. It also appears in the Klausmeier model of vegetation and precipitation (section \ref{sec:klaus}) and the harvested predator-prey system (section \ref{sec:predprey}). In the Klausmeier model, an increase in the bifurcation value of the rainfall parameter $\p$ as return time $\tau$ increases means that that higher average precipitation rates are needed to maintain a stable vegetated equilibrium as the time between precipitation events becomes longer. In the predator-prey model, sustainable harvesting of the predator population requires lower harvesting rates when harvests are discretized, more spread out, and more intense, because the transcritical bifurcation that represents predator collapse occurs at a lower value of the  harvesting parameter $\p$ as the recovery period $\tau$ increases.
As noted by \cite{meyer2018quantifying} and \cite{zeeman2018resilience}, these patterns have real implications for a system's resilience to disturbance. First, we may underestimate a system's vulnerability to regime changes if we smooth out the discrete quality of disturbances in our models. Second, changes to the timing of real-world disturbances---while maintaining the same average disturbance rate---can drive tipping. Importantly, the magnitude of repeated kicks needed to drive tipping tends to be lower than the distance to a threshold in state space, due to the accumulation of disturbance and incomplete recovery. While the details will vary by system, our bifurcation continuation results may aid analysis of flow-kick bifurcation thresholds by providing a seed for numerical continuation starting near $\tau=0$. 

In addition, proportionally increasing kick magnitudes and return periods can introduce entirely new dynamics not observed in the continuous analog system. 
We observed this when period-doubling bifurcations arose in the flow-kick models of the predator-prey example (Section \ref{sec:predprey}). Although these occurred in an infeasible parameter space, we suspect period doubling or chaos may also appear in feasible regions at resonance bifurcations along the Neimark-Sacker curve. In addition,
in a discrete model of fire disturbances to a tree-grass savanna, Hoyer-Leitzel and Iams found curves of bifurcations in disturbance space that are isolated away from the origin---bifurcations that occur in the flow-kick systems but not in their continuous analog \cite{hoyer2021impulsive}. It should be noted that the kick structure in \cite{hoyer2021impulsive} differs from Definition \ref{def:fkr}, with the disturbance rate dependent on $\tau$ as well as $x$. Broadly, these examples alert  us that accounting for the discrete nature of disturbance could reveal bifurcation structures absent from an ODE model. 
A rich area for future investigation includes studying the dynamic structures that flow-kick models gain relative to their ODE counterparts.

\bibliographystyle{siamplain}
\bibliography{ctslimit.bib}

\pagebreak

\appendix
\section{Integral Form of Remainder}\label{sec:intform}
The following integral form of the remainder of a Taylor series is standard (e.g.
\cite{folland1990remainder}):

\begin{proposition}\label{prop:remain} If $f$ is of class $C^{k+1}$ on an open interval $I\subset\R$ with a base point $a\in I$ and the variable $t\in I$, then 
\begin{equation*}
    f(t)=\sum_{j=0}^{k} \frac{f^{(j)}(a)}{j!}(t-a)^j+R_{k,a}(t)\
\end{equation*}
where $\displaystyle R_{k,a}(t) =\frac{(t-a)^{k+1}}{k!}\int_0^1 (1-v)^k f^{(k+1)}(a+v(t-a))dv$.
\end{proposition}

When we expand $\varphi(\tau,x)$ in $\tau$ in Proposition \ref{prop:fixedpt_n} and thereafter, we use Proposition \ref{prop:remain} with $k=1$, $a = 0$, $t = \tau$, and $f(\cdot)=\varphi(\cdot,x)$.

\newpage

\section{Further details on the proofs of bifurcation results}\label{sec:bif_app}

\subsection{Details from Theorem \ref{thm:s-n}}\label{sec:bif_sn}
\noindent \textbf{Confirming that} $\mathbf{\tfrac{\partial r_1}{\partial \p}(x_*,y_*,\p*) \neq 0}$. \\
We have by our assumption on the vector field in Table \ref{table:SNnew} that
\begin{align}
    0&\not=w_*\frac{\partial}{\partial \p}[\tilde f+ \tilde r](\tilde x_*,\p*)\\
    &=w_*\frac{\partial \tilde r}{\partial \p}(\tilde x_*,\p*) & \text{(since $\partial \tilde f/\partial \p=0$)}\\
    &=w_*\frac{\partial}{\partial \p}[r+D_{\tilde x}\tilde r(\tilde x_*,\p_*)\tilde x]& \text{(definition of $r$)}\\
    &=w_*\frac{\partial r}{\partial \p}\\
    &=\frac{\partial r_1}{\partial \p}(x_*,y_*,\p*) & \text{(since $w_*=[1 \ 0 \cdots 0]$).}   
\end{align}

\noindent \textbf{Confirming that} $\mathbf{\tfrac{\partial F_3}{\partial x}(0,x_*,y_*,\p*) \neq 0}$. \\
First recall that $F_3 = \text{det}(D_{(x,y)}M)$ where 
\begin{equation}\label{eq:F3_start}
M = \begin{bmatrix} f_1(x,y)+\mathcal{O}(\tau)+r_1( \varphi_\tau(x,y),\p) \\
By + f_{2\text{-}n}(x,y)+\mathcal{O}(\tau) + r_{2\text{-}n}(\partial \varphi(x,y),\p)\end{bmatrix}.
\end{equation}
We compute
\begin{align}\label{eq:rearr}
    D_{(x,y)}M &=     \left[ \begin{array}{c|ccc}
    \frac{\partial f_1}{\partial x} + \frac{\partial r_1}{\partial x}+\frac{\partial}{\partial x}\mathcal{O}(\tau) 
    & \cdots & \frac{\partial f_1}{\partial y_j}+\frac{\partial r_1}{\partial y_j}+\frac{\partial}{\partial y_j}\mathcal{O}(\tau)  & \cdots  
    \\
    \hline
    \vdots & \ddots & \vdots & \iddots \\
     \frac{\partial f_{i+1}}{\partial x}+\frac{\partial r_{i+1}}{\partial x}+\frac{\partial}{\partial x}\mathcal{O}(\tau) & \cdots &B_{ij}\!+\!\frac{\partial f_{i+1}}{\partial y_j}\!+\frac{\partial r_{i+1}}{\partial y_j}\!+\!\frac{\partial}{\partial y_j}\mathcal{O}(\tau) & \cdots\\
    \vdots & \iddots  &\vdots &\ddots \\
    \end{array}\right], 
\end{align}
noting that $\partial r_k/\partial x=\nabla_u r_k\cdot(\partial \varphi_\tau/\partial x)$ and similarly for $\partial r_k/\partial y_j$.
Performing a cofactor expansion across the first row, we get
{\allowdisplaybreaks
\begin{align}
\begin{split}
\det(D_{(x,y)}M)&=\Bigg[ \left(\frac{\partial f_1}{\partial x} + \frac{\partial r_1}{\partial x} +\frac{\partial}{\partial x}\mathcal{O}(\tau) \right)\det\left(B+D_y\left[f_{2\text{-}n}+r_{2\text{-}n} \right] + D_y\mathcal{O}(\tau) \right) 
\end{split}\\
\begin{split}
    & \hspace{1cm} +\sum\limits_{j=1}^{n-1}(-1)^j\left(\frac{\partial f_1}{\partial y_j}+ \frac{\partial r_1}{\partial y_j} +\frac{\partial}{\partial y_j}\mathcal{O}(\tau) \right) \det(M_{1,j+1} )\Bigg]\label{eq:det_end}
\end{split}\\
    &= F_3(\tau,x,y,\p)\label{eq:F3_end}
\end{align}
}
where $M_{1,j+1}$ is the $(1,j+1)$ minor matrix of $D_{(x,y)}M$. Condensing notation for $F_3$, let
\begin{align}
    \frac{\partial f_1}{\partial x} + \frac{\partial r_1}{\partial x} +\frac{\partial}{\partial x}\mathcal{O}(\tau)
    &=\alpha(\tau,x,y,\p),\label{eq:alpha}\\
    \det\left(B+D_y\left[f_{2\text{-}n}+r_{2\text{-}n} \right] + D_y\mathcal{O}(\tau) \right)
    &=\beta(\tau,x,y,\p),\label{eq:beta}\\
    \text{and }\quad\quad\frac{\partial f_1}{\partial y_j}+ \frac{\partial r_1}{\partial y_j} +\frac{\partial}{\partial y_j}\mathcal{O}(\tau)
    &=\gamma_j(\tau,x,y,\p),\label{eq:gamma}
\end{align}
so that
\begin{equation}
    F_3=\alpha \, \beta+\sum\limits_{j=1}^{n-1}(-1)^j\gamma_j\det(M_{1,j+1}).
\end{equation}
Now we compute the derivative of $F_3$ with respect to $x$:
\begin{equation}
    \frac{\partial F_3}{\partial x} = \frac{\partial \alpha} {\partial x} \beta + \alpha \frac{\partial \beta}{\partial x} + 
    \sum_{j=1}^{n-1}(-1)^j \frac{\partial \gamma_j}{\partial x} \det(M_{1,j+1}) + \sum_{j=1}^{n-1} (-1)^j \gamma_j \frac{\partial }{\partial x}\left[\det(M_{1,j+1})\right].
\end{equation}

The term $\partial\alpha/\partial x$ is nonzero at $(0,x_*,y_*,\p_*)$ due to the quadratic dominance assumption on the vector field $f+r$ in Table \ref{table:SNnew}. In particular, when $\tau=0$ we have $ \varphi_0= \ $id, and
\begin{align}
    \frac{\partial}{\partial x}\left[\alpha(0,x,y,\p)\right]
    &=\frac{\partial^2 f_1}{\partial x^2} + \frac{\partial^2 r_1}{\partial x^2}\\
    &=
    \frac{\partial^2 f_1}{\partial x^2}+\frac{\partial}{\partial x}\left[\nabla_u r_1\cdot\frac{\partial \varphi_0}{\partial x}\right]\\
    &=
    \frac{\partial^2 f_1}{\partial x^2}+\frac{\partial}{\partial x}\left[\nabla_u r_1\right]\cdot\frac{\partial \varphi_0}{\partial x}+
    \nabla_u r_1 \cdot \frac{\partial^2 \varphi_0}{\partial x^2}\\
    &=
    \frac{\partial^2 f_1}{\partial x^2}+\left[\cdots \nabla_u \frac{\partial r_1}{\partial u_k} \cdot \frac{\partial  \varphi_0}{\partial x} \cdots \right]\cdot\frac{\partial \varphi_0}{\partial x}+
    \nabla_u r_1 \cdot \frac{\partial^2 \varphi_0}{\partial x^2}.
\end{align}
Using the facts that $\partial \varphi_0/\partial x=[1 \ 0 \cdots \ 0]^T$, $\partial^2 \varphi_0/\partial x^2=[0 \ 0 \cdots \ 0]^T$, and $u_1\equiv x$, we get
\begin{align}
    \frac{\partial}{\partial x}\left[\alpha(0,x,y,\p)\right]
    &=\frac{\partial^2 f_1}{\partial x^2}+\frac{\partial^2 r_1}{\partial x^2}.\label{eq:alpha_end}
\end{align}
On the other hand, from Table \ref{table:SNnew} we have
\begin{align}
    w_*(D_x^2[\tilde f+ \tilde r](x_*,\p_*)(v_*,v_*))  &\not=0\label{eq:alpha_start} \\
    \implies
    w_*\sum\limits_{i,j=1}^n\frac{\partial^2(\tilde f+ \tilde r)}{\partial \widetilde x_i \partial \widetilde x_j}\Bigg|_{\widetilde x_*,\p_*}v_{*i}v_{*j}&\not=0. \\
    \implies
    w_*\sum\limits_{i,j=1}^n\frac{\partial^2}{\partial \widetilde x_i \partial \widetilde x_j}\left[f+D_{\tilde x}\tilde f(\tilde x_*,\p_*)\tilde x+r+D_{\tilde x}\tilde r(\tilde x_*,\p_*)\tilde x \right]\Bigg|_{\widetilde x_*,\p_*}v_{*i}v_{*j}&\not=0\\
    \implies
    w_*\sum\limits_{i,j=1}^n\frac{\partial^2}{\partial \widetilde x_i \partial \widetilde x_j}\Big[f+r\Big]\Bigg|_{\widetilde x_*,\p_*}v_{*i}v_{*j}&\not=0.\label{eq:w*v*}
\end{align}
Since $w_*=[1 \ 0 \ \cdots \ 0 ]$ and $v_*=[1 \ 0 \ \cdots \ 0 ]^T$, \eqref{eq:w*v*} reduces to
\begin{equation}\label{eq:B1_red}
    \frac{\partial^2(f_1+r_1)}{\partial \tilde x_1^2}\Bigg|_{x_*,y_*,\p_*}\not=0. 
\end{equation}
Comparing equations \eqref{eq:alpha_start}-\eqref{eq:alpha_end} with equation \eqref{eq:B1_red} and recalling $\tilde x_1=x$, we conclude that when $(\tau,x,y,\p)= (0,x_*,y_*,\p_*)$ the expression $\partial \alpha/\partial x$ is nonzero.

The term $\beta$ reduces to the nonzero quantity $\det(B)$ at $(0,x_*,y_*,\p_*)$. Furthermore, $\alpha$ and $\gamma$ evaluate to $0$ at $(0,x_*,y_*,\p_*)$ (see equation \eqref{eq:sn3}). Lastly, $\det(M_{i,j+1})$ is zero at $(0,x_*,y_*,\p_*)$ because the first column of each minor matrix $M_{1,j+1}$ consists entirely of entries $\frac{\partial f_{i+1}}{\partial x}+\frac{\partial r_{i+1}}{\partial x}+\frac{\partial}{\partial x}\mathcal{O}(\tau)$ that vanish at $(0,x_*,y_*,\p_*)$. Together this implies that $\frac{\partial F_3}{\partial x}(0,x_*,y_*,\p_*)\not=0$.\\

\noindent\textbf{(iii) Non-degenerate eigenvalues continue}\\
We will show that if $n-1$ eigenvalues of $D_x[f+r](x_*,\p_*)$ have nonzero real parts, then $n-1$ eigenvalues of $D_{(x,y)}\Phi_\tau(x_\tau,y_\tau,\p_\tau)$ do not have modulus 1 along the branch of flow-kick fixed points $(x_\tau,\p_\tau)$.

As shown in equations \eqref{eq:DxyPhitau}-\eqref{eq:M},
\begin{equation}
    D_{(x,y)}\Phi_\tau=I+\tau D_{(x,y)}\begin{bmatrix} f_1(x,y)+\mathcal{O}(\tau)+r_1( \varphi_\tau(x,y),\p) \\
By + f_{2\text{-}n}(x,y)+\mathcal{O}(\tau) + r_{2\text{-}n}( \varphi_\tau(x,y),\p)\end{bmatrix}
\end{equation}
and so
\begin{align}\label{eq:rearr-ed}
    D_{(x,y)}\Phi_\tau\Bigg|_{\tiny\begin{array}{l}x_\tau\\y_\tau \\ \p_\tau \end{array}} &=I+\tau
    \underbrace{\left[ \begin{array}{c|ccc}
    \frac{\partial f_1}{\partial x} + \frac{\partial r_1}{\partial x} +\frac{\partial}{\partial x}\mathcal{O}(\tau) 
    & \cdots & \frac{\partial f_1}{\partial y_j}+\frac{\partial r_1}{\partial y_j}+\frac{\partial}{\partial y_j}\mathcal{O}(\tau)  & \cdots  
    \\
    \hline
    \vdots & \ddots & \vdots & \iddots \\
     \frac{\partial f_{i+1}}{\partial x}+\frac{\partial r_{i+1}}{\partial x}+\frac{\partial}{\partial x}\mathcal{O}(\tau) & \cdots &B_{ij}+\frac{\partial f_{i+1}}{\partial y_j}+\frac{\partial r_{i+1}}{\partial y_j}+\frac{\partial}{\partial y_j}\mathcal{O}(\tau) & \cdots\\
    \vdots & \iddots  &\vdots &\ddots \\
\end{array}\right]}_{A(\tau,x_\tau,y_\tau,\p_\tau)}. 
\end{align}
Including multiplicities, denote the eigenvalues of $ D_{(x,y)}\Phi_\tau(x_\tau,y_\tau,\p_\tau)$ as $\ell_1(\tau),\cdots,\ell_n(\tau)$ and the eigenvalues of  $A(\tau,x_\tau,y_\tau,\p_\tau)$ as $m_1(\tau), \cdots m_n(\tau)$. We can see that 
\begin{align}\label{eq:rearr-0}
    A(0,x_*,y_*,\p_*)&=
   \left[ \begin{array}{c|ccc}
    0
    & \cdots & 0 & \cdots  
    \\
    \hline
    \vdots &  \\
     0 &   &B  &  \\
    \vdots &  \\
\end{array}\right] 
\end{align}
as follows. At $\tau=0$ we have $\varphi_0(x,y)=(x,y)$ and $(x_0,y_0)=(x_*,y_*)$, so
\begin{align}
    \frac{\partial r_k}{\partial \tilde x_\ell}( \varphi(x_0,y_0),\p_0)&=\frac{\partial r_k}{\partial \tilde x_\ell}(x_*,y_*,\p_*)\\
    &=\nabla_u r_k(x_*,y_*,\p_*) \frac{\partial  \varphi_0}{\partial \tilde x_\ell}(x_*,y_*) \\
    &=\left[\cdots \partial r_k/\partial u_m (x_*,y_*,\p_*)\cdots \right][0 \cdots 0 \underbrace{ \ 1 \ }_{\ell^\text{th}\text{ entry}} 0 \cdots 0]^T\\
    &=\partial r_k/\partial u_\ell(x_*,y_*,\p_*)\\
    &\equiv \partial r_k/\partial \tilde x_\ell(x_*,y_*,\p_*).
\end{align}
It follows that
\begin{align}
    \frac{\partial}{\partial \tilde x_\ell}\Big[f_k(\tilde x)+r_k( \varphi_0(\tilde x),\p_0)\Big](x_*,y_*,\p_*)
    &=\frac{\partial}{\partial \tilde x_\ell}\Big[f_k(\tilde x)+r_k(\tilde x,\p_0)\Big](x_*,y_*,\p_*)\\
    &=0
\end{align}
by \eqref{eq:sn3}-\eqref{eq:sn4}. Substituting into the matrix $A(\tau,x_\tau,y_\tau,\p_\tau)$ in equation \eqref{eq:rearr-ed} and noting $\mathcal{O}(\tau)$ vanishes at $\tau=0$ finishes the confirmation of \eqref{eq:rearr-0}.

By assumption, $\text{Re}(\text{spec} \ B)\not=0$, so $A(0,x_*,y_*,\p_*)$ has eigenvalues $m_1(0)=0$ and $\text{Re}(m_i(0))\not=0$ for $i=2,\cdots,n$. Since the eigenvalues of $A$ are continuous with respect to its entries \cite{stewart1990matrix}, which are continuous functions of $\tau$, we retain the condition $\text{Re}(m_i(\tau))\not=0$ for sufficiently small $\tau$ and $i=2,\cdots,n$. It follows from the fact that $\ell_i(\tau)=1+\tau m_i(\tau)$ and the argument used in the proof of Theorem \ref{thm:stability_n} that $|\ell_i(\tau)|\not=1$, again for sufficiently small $\tau$ and $i=2,\cdots,n$.\\

\noindent \textbf{(iv) Transversality}\\
We use the notation $x\in\R^n$ rather than the notation $\tilde x=(x,y)\in\R^n$.

Transversality of the flow-kick map with respect to parameter $\p$ requires that 
\begin{equation}
    w_\tau \frac{\partial \Phi_\tau}{\partial \p} (x_\tau,\p_\tau)\not=0, 
\end{equation}
which is equivalent to the condition
\begin{equation}
\tau w_\tau \left(\frac{\partial r}{\partial \p}(\varphi_\tau(x_\tau),\p_\tau)\right)\not=0.
\end{equation}
We have $\tau>0$ for flow-kick maps $\Phi_\tau$. 
To see that $w_\tau \left(\frac{\partial r}{\partial \p}(\varphi_\tau(x_\tau),\p_\tau)\right)$ is nonzero for $\tau>0$ sufficiently small, first note that $r$ is $C^2$ and $(\varphi_\tau(x_\tau),\p_\tau)$ is continuous in $\tau$, so $\partial r/\partial\p(\varphi_\tau(x_\tau),\p_\tau)$ is continuous in $\tau$. On the one hand, since $w_\tau$ is a left eigenvector of $D_{x}\Phi_\tau=I+\tau A(\tau,x_\tau,\p_\tau)$ (see equation \eqref{eq:rearr-ed}) with eigenvalue $1$, we have
\begin{align}
w_\tau&\in \ker((I+\tau A(\tau,x_\tau,\p_\tau))^T-1I)\\
&=\ker(A(\tau,x_\tau,\p_\tau)^T) &\text{for } \tau>0.
\end{align}
On the other hand, the left eigenvalue $w_*$ of
\begin{equation}
D_{x}\Big[f(x)+r(x,\p)\Big]\Bigg|_{(x_*,\p_*)}=
   \left[ \begin{array}{c|ccc}
    0
    & \cdots & 0 & \cdots  
    \\
    \hline
    \vdots &  \\
     0 &   &B  &  \\
    \vdots &  \\
\end{array}\right] 
=A(0,x_*,\p_*)
\end{equation}
associated with eigenvalue $0$ satisfies
\begin{align}
    w_*&\in\ker(A(0,x_*,y_*,\p_*)^T-0I)\\
    &=\ker(A(0,x_*,\p_*)^T).
\end{align}
The matrix $A(\tau,x_\tau,\p_\tau)^T$ is continuous in $\tau$ and evaluates to $A(0,x_*,\p_*)^T$ at $\tau=0$. 
It follows from continuity of nullspaces with respect to matrix perturbations 
that $w_\tau$ varies continuously in $\tau$ from $w_*$ at $\tau=0$, up to nonzero scaling. Up to nonzero scaling, the product $w_\tau \left(\frac{\partial r}{\partial \p}(x_\tau,\p_\tau)\right)$ is also continuous in $\tau$ and evaluates to the quantity
$
    w_*\left(\frac{\partial r}{\partial \p}(x_*,\p_*)\right)
$
at $\tau=0$, which is nonzero by the transversality assumption on the ODE system (see Table \ref{table:SNnew}) and the fact $\partial f/\partial\p=0$ . Hence $\tau w_\tau \left(\frac{\partial r}{\partial \p}(x_\tau,\p_\tau)\right)$ is nonzero for $\tau>0$ sufficiently small.\\

\noindent\textbf{(v) Quadratic dominance}\\
Again, use original coordinates $x\in\R^n$. We will show that  $w_\tau (D_x^2\Phi_\tau(x_\tau,\p_\tau)(v_\tau,v_\tau))\not=0 $
for $\tau$ sufficiently small. Using the coordinate definition of the quadratic tensor product and the Taylor expansion for $\Phi_\tau(x,\p) = x+\tau f(x)+\tau r(\varphi_\tau(x),\p)+\mathcal{O}(\tau^2)$, we obtain
\allowdisplaybreaks{
\begin{align}
    &\phantom{=} w_\tau (D_x^2\Phi_\tau(x_\tau,\p_\tau)(v_\tau,v_\tau) ) = w_\tau \sum_{i,j=1}^n \left.\frac{\partial^2 \Phi_\tau}{\partial x_i \partial x_j}\right|_{x_\tau,\p_\tau} v_{\tau i} v_{\tau j} \\
    & =  w_\tau \sum_{i,j=1}^n \left.\left(\frac{\partial x}{\partial x_i \partial x_j} + \tau\frac{\partial^2 f}{\partial x_i \partial x_j} + \tau\frac{\partial^2}{\partial x_i \partial x_j}\Big[ r(\varphi_\tau(x),\lambda) \Big]+ \mathcal{O}(\tau^2)\right)\right|_{x_\tau,\p_\tau} v_{\tau i} v_{\tau j} \\
    &  = \tau w_\tau \sum_{i,j=1}^n 
    \left.\left( \frac{\partial^2 f}{\partial x_i \partial x_j} + \frac{\partial}{\partial x_i}\Big[D_u r(\varphi_\tau(x),\p) \Big]\frac{\partial \varphi_\tau}{\partial x_j}  + D_u r(\varphi_\tau(x),\p) \frac{\partial^2 \varphi_\tau}{\partial x_i\partial x_j}+ \mathcal{O}(\tau) \right)
     \right|_{x_\tau,\p_\tau} \!\!\! v_{\tau i} v_{\tau j} 
    \\
    &=
   \tau w_\tau \sum_{i,j=1}^n 
    \left.\left( \frac{\partial^2 f}{\partial x_i \partial x_j} + \begin{bmatrix} & : & \\
    \cdot\cdot & \nabla_u\tfrac{\partial r_k}{\partial u_\ell}\cdot\tfrac{\partial\varphi_\tau}{\partial x_i} &\cdot\cdot\\
    & : & \end{bmatrix}\frac{\partial \varphi_\tau}{\partial x_j}  + D_u r(\varphi_\tau(x),\p) \frac{\partial^2 \varphi_\tau}{\partial x_i\partial x_j}+ \mathcal{O}(\tau) \right)
     \right|_{x_\tau,\p_\tau} \!\!\!\!\!\!\!\!\!\! v_{\tau i} v_{\tau j} 
    \\
    & = \tau G(\tau,x_\tau,\lambda_\tau).
    \end{align}
    }
We will show that when $\tau =0$, $G(0,x_0,\lambda_0)=G(0,x_*,\p_*) \neq 0$. 

When $\tau = 0$, eigenvectors evaluate to $w_\tau = w_* = \begin{bmatrix} 1 &  0& \cdots & 0 \end{bmatrix}$ and $v_\tau = v_* = \begin{bmatrix} 1 & 0 & \cdots & 0\end{bmatrix}^T$. Furthermore, $\varphi_\tau(\cdot) = \varphi_0(\cdot) = \text{id}(\cdot)$, making  $\tfrac{\partial \varphi_0}{\partial x_j}$  the $j$th standard basis vector $e_j$ and $\tfrac{\partial^2 \varphi_0}{\partial x_i \partial x_j} = $ for all $i,j$.
It follows that
\begin{align}
  G(0,x_*,\lambda_*)  & = w_* \left.\left( \frac{\partial^2 f}{\partial x_1^2} + \begin{bmatrix} & : & \\
    \cdot\cdot & \nabla_u\tfrac{\partial r_k}{\partial u_\ell}\cdot e_1 &\cdot\cdot\\
    & : & \end{bmatrix}\begin{bmatrix}1 \\ 0 \\ : \\ 0\end{bmatrix} + D_u r(\varphi_\tau (x),\p) \, 0 + 0\right)\right|_{x_*,\p_*}  \\
   & \equiv w_* \left.\left(\frac{\partial ^2 f}{\partial x_1^2} + \frac{\partial^2 r}{\partial x_1^2} \right)\right|_{x_*,\p_*}\label{eq:quad_dom}
\end{align}
and one can check that \eqref{eq:quad_dom} is equivalent to $w_*(D_x^2[f+r](x_*,\p_*)(v_*,v_*))$, a nonzero quantity by assumption.

Since $w_\tau, v_\tau, f, r, \varphi$ are all $C^2$ functions, $G$ is a continuous function of $\tau$. Since when $\tau =0$, $G(0,x_*,\p_*)\neq 0$,  there exists a  neighborhood of $\tau =0$ in which $\tau G(\tau,x_\tau,\p_\tau)$ is also nonzero.

\vspace*{0.5cm}

\subsection{Details from Theorem \ref{thm:tc}}\label{sec:bif_tc}
\noindent \textbf{Confirming $\mathbf{\frac{\partial^2 r_1}{\partial \p \partial x}(x_*,y_*,\p_*)\not=0}$}\\

From the mixed partials condition on the continuous system in Table \ref{table:TC}, we have
\begin{align}
    0&\not=w_*\left(D_{\tilde x}\frac{\partial}{\partial \p}\left[\tilde f+ \tilde r\right](\tilde x_*,\p_*)v_*\right)\\
    &=w_*\left(D_{\tilde x}\frac{\partial}{\partial \p}\left[r + D_{\tilde x}(\tilde x_*,\p_*)\tilde x\right](\tilde x_*,\p_*)v_*\right)\\
    &=w_*\left(D_{\tilde x}\frac{\partial r}{\partial \p}(\tilde x_*,\p_*)v_*\right)\\
    &=\frac{\partial^2 r_1}{\partial x\partial \p}(\tilde x_*,\p_*)\\
    &=\frac{\partial^2 r_1}{\partial \p\partial x}( x_*,y_*,\p_*)
\end{align}
as desired.\\

\noindent \textbf{(v) Mixed partials}\\
To confirm the mixed partials condition (v), we work in original coordinates ($x\in\R^n$). Consider
\begin{align}
    w_\tau\left(D_x\frac{\partial \Phi_\tau}{ \partial \p}(x_*,\p_\tau)v_\tau\right)
    &=\tau w_\tau\left(D_x\frac{\partial r}{\partial \p}(x_*,\p_\tau)v_\tau\right).
\end{align}
Clearly $\tau>0$ for flow-kick maps. From the second factor define
\begin{equation}
    h(\tau)=w_\tau\left(D_x\frac{\partial r}{\partial \p}(x_*,\p_\tau)v_\tau\right).
\end{equation}
From similar arguments to those used in Theorem \ref{thm:s-n} we get that $h$ is a continuous function of $\tau$, and
\begin{align}
    h(0)&=w_*\left(D_x\frac{\partial r}{\partial \p}(x_*,\p_*)v_*\right).
\end{align}

From the mixed partials condition (v) for the vector field in Table \ref{table:TC}, we have
\begin{align}
    0&\not=w_*\left(D_x\frac{\partial }{\partial \p}[f+r](x_*,\p_*)v_*\right).
\end{align}

Since $f$ does not depend on $\p$, $h(0)$ is equivalent to the mixed partials expression (v) for the vector field (see Table \ref{table:TC}); namely,
\begin{align}
    h(0)&=w_*\left(D_x\frac{\partial }{\partial \p}[f+r](x_*,\p_*)v_*\right)\not=0.
\end{align}
The function $h(\tau)$ remains nonzero for sufficiently small $\tau$, and so does $\tau w_\tau\left(D_x\frac{\partial \Phi_\tau}{ \partial \p}(x_*,\p_\tau)v_\tau\right)$.

\newpage 
\section{Conditions for bifurcations in the Predator Prey system\label{sec:predprey_bifconditions}} 
For the predator-prey flow-kick map given by equation \eqref{eq:predprey_flowkick}, we derive conditions for bifurcations at the point (4,0).

\textit{Derivation:}
    Let $\Phi_{\tau}(x,y) = 
    \begin{bmatrix}
        [\varphi_\tau(x,y)]_1 \\ [\varphi_\tau(x,y)]_2-\tau\p [\varphi_\tau(x,y)]_2
    \end{bmatrix} \equiv
    \begin{bmatrix}
        U(x,y) \\ V(x,y)
    \end{bmatrix}$.
    We have used subscripts of $1$ and $2$ on the flow function $\varphi$ to denote its $x$ and $y$ components, respectively. For all values of the disturbance parameters $\tau$ and $\p$, (4,0) is a fixed point for the map; that is,
    \begin{equation}
      \left\{  \begin{array}{ll}
          4   & = U(4,0)  \\
          0   & = V(4,0).
        \end{array} \right.
    \end{equation}

To find bifurcations, we will look for conditions at which an eigenvalue of the Jacobian 
\begin{equation} \left. \begin{bmatrix} \partial U/\partial x & \partial U/\partial y \\
\partial V /\partial x & \partial V/\partial y \end{bmatrix} \right|_{(x,y)=(4,0)} \end{equation}
of $\Phi$ has absolute value equal to 1. 
We first note that the set $y = 0$ is invariant for both the flow and the map. Thus $\partial V/\partial x |_{y=0} = 0$ and the Jacobian is an upper triangular matrix with real valued eigenvalues equal to $\partial U/\partial x |_{(x,y)=(4,0)}$ and $\partial V/\partial y |_{(x,y)=(4,0)}$. Since $y=0$ is invariant, we can calculate $U(x,0)$ exactly and find \begin{equation} U(x,0) = \dfrac{4x}{x+(4-x)e^{-\tau}}
\end{equation}
and $\partial U/\partial x |_{(4,0)} = e^{-\tau}.$ For $\tau > 0$, this eigenvalue is always less than 1. Thus bifurcations depend only on the second eigenvalue.

We estimate the second eigenvalue by approximating $V(4,y)$.
We fix the value of $x$ at $4$ in the second equation of \eqref{eq:predprey_undisturbed}, giving \begin{equation}\label{eq:dydt_x=4}
    \begin{aligned}
        \frac{dy}{dt} &= -0.5y+y(1-e^{-2}) \\
        &= (1/2-e^{-2})y.
    \end{aligned}
    \end{equation}
Integrating \eqref{eq:dydt_x=4} yields $y(t)=y_0e^{a\tau}$ where $a=(1/2-e^{-2})$. We use  $ye^{a\tau}$ to approximate $[\varphi_\tau(4,y)]_2$, noting that the approximation is exact at $y=0$. Then
    \begin{equation}\label{eq:appCv}
        V(4,y)\approx ye^{a\tau}-\tau\p ye^{a\tau}=(1-\tau\p)ye^{a\tau}
    \end{equation}
and  $\partial V/\partial y (4,y) \approx (1-\tau \p)e^{a\tau}$. 

With the approximation \eqref{eq:appCv}, we see the Jacobian at $(4,0)$ has an eigenvalue of $1$ when 
\begin{equation}
    (1-\tau\p)e^{a\tau} =1 \label{eq:tc_curve2}
\end{equation}
This indicates that there could be a bifurcation along this curve. We can verify that the curve meets necessary conditions for a transcritical bifurcation when the curve $y-V(4,y) = 0$ has $y=0$ as a double root. We rewrite $y-V(4,y)=0$ as
    \begin{align}
        y-(1-\tau\p)ye^{a\tau}&=0\\
        \iff y(1-(1-\tau\p)e^{a\tau})&=0.\label{eq:appCexact}
    \end{align}
    Now at $y=0$ equation \eqref{eq:appCexact} is equivalent to the fixed point condition $y-V(4,y)=0$. Requiring that this equation has a double root at $y=0$ amounts to requiring
    \begin{align}
        1-(1-\tau\p)e^{a\tau}&=0\\
        \iff \tau\p=1-e^{-a\tau}.\label{eq:tc_curve}
    \end{align}
Noting that the curve of possible bifurcations given by \eqref{eq:tc_curve2} is equivalent to \eqref{eq:tc_curve},
we expect to find a transcritical bifurcation in $\tau\lambda$-space along this curve.

Additionally, we find the system undergoes a period doubling bifurcation when the eigenvalue is $-1$, which occurs where 
\begin{equation} \partial V/\partial y (4,0)=(1-\tau\lambda)e^{a\tau} =  -1. \end{equation}
Solving for $\tau\lambda$, we get 
\begin{equation}
    \tau\lambda = 1+e^{-a\tau} >1 \label{eq:period_doubling}
\end{equation}
which represents a harvesting percentage greater than 100\% of the predator population. Therefore the period doubling bifurcation occurs in a nonphysical region of parameter space.

\end{document}